\documentstyle[twoside,11pt]{article}

\def\R#1{\par\hangindent\parindent\indent\llap{#1\enspace}\ignorespaces}
\begin{document}
\title{
{\bf Shape of Alexandrov spaces
\\
with positive Ricci curvature}
\author { Zisheng Hu   \ \ \ \ \ \ \   Le Yin}
\date {}}
\maketitle

\begin{center}
\begin{minipage}{135mm}
{\bf Abstract:} Under the definition of Ricci curvature bounded below for Alexandrov spaces introduced by Zhang-Zhu,
we extend a result by Colding that an $n-$dimensional manifold with Ricci curvature greater or equal to $n-1$
and volume close to that of the unit $n-$sphere is close (in the Gromov-Hausdorff distance) to the sphere,
from the case of Riemannian manifolds to the case of Alexandrov spaces,
with an additional assumption, roughly speaking, that the rough volume of the set of "short" cut points is small,
following the basic idea in the Riemannian case with necessary modifications
because of the only almost everywhere second differentiability of distance functions.
\\
%{\bf Key words:}
%\\
%{\bf MR(2000) Subject Classification:}
%\\
\end{minipage}
\end{center}

\bigbreak

{\section {Introduction }}

Let $\omega_{n}$ be the volume of the round $n-$sphere $S^{n}$ with sectional curvature one, in [C], Colding shows the following
\\

\noindent
{\bf Theorem 1.1 [C]}
\textit{
Given an integer $n\geq 2$ and an $\epsilon>0$, there exists a $\delta(n, \epsilon)>0$ such that if $M^{n}$ is an $n-$dimentional manifold with
$Ric(M)\geq n-1$ and $Vol(M)\geq\omega_{n}-\delta$, then the Gromov-Hausdorff distance between $M$ and $S^{n}$ is at most $\epsilon$.
}
\\

As the concepts of curvature in Riemannian geometry have been generalized to singular spaces,
it would be interesting to extend the above result to the case of singular spaces.

Below we briefly recall the generalizations of lower bound of sectional curvature and Ricci curvature in Riemannian geometry.
As for sectional curvature, the generalization is Alexandrov spaces with curvature bounded below,
see the seminal paper [BGP] and the 10th chapter in the text book [BBI];
As for Ricci curvature, there are several notions to generalize lower Ricci curvature bound to metric measure spaces,
such as Lott-Sturm-Villani's curvature-dimension condition $CD(k, n)$ [LV, S1, S2], Ohta-Sturm's $MCP$ condition[O, S2] and Ambrosio-Gigli-Savar\'{e}'s
Riemannian curvature dimension condition $RCD(k, n)$ [AGS, EKS].

In this paper, we are concerned with lower Ricci curvature bound for Alexandrov spaces as a particular case of metric measure spaces.
In [KS1, KS2], Kuwae-Shioya introduced an infinitesimal Bishop-Gromov condition $BG(k)$ for Alexandrov spaces, and established a topological version of Cheeger-Gromoll splitting theorem on Alexandrov spaces under $BG(k)$;
In [ZZ1], Zhang-Zhu  introduced a definition of Ricci curvature bounded below for Alexandrov spaces based on Petrunin's second variation of arc length [Pet1],
they showed that the new definition implies the curvature-dimension condition
and there hold Cheeger-Gromoll splitting theorem (see also [G] on general metric measure spaces) and maximal diameter theorem on Alexandrov spaces under this Ricci condition.
Furthermore, they established a Bochner type formula on Alexandrov spaces and many important results under this Ricci condition, see [ZZ2, ZZ3, QZZ, ZZ4].

See [ZZ2] for a detailed discussion of the relations among various generalizations of Ricci curvature on Alexandrov spaces.

To state the main result, denote

$M$ --- an $n-$dimensional Alexandrov space with curvature bounded below by some $k_{0}(<0)$ and without boundary;

$C_{p, \hat{\theta}}$ --- the set of points $x\in M$ such that minimal geodesic $\overline{px}$ can not be extended beyond $x$
  and $|px|\leq \pi-\hat{\theta}$, \ for any $p\in M$, $\hat{\theta}\rightarrow 0^{+}$;

$\beta_{X}(a)$ --- the largest possible number of points $x_{i}\in X$ that are at least $a$ pairwise distant from each other,
\ for any metric space $X$ and $a>0$;

$Vr_{d}(X):= \limsup_{a\rightarrow 0^{+}}a^{d}\beta_{X}(a)$ --- the rough $d-$dimensional volume of a bounded set $X$ in a metric space,
for any $d\geq 0$.
\\

Under the definition of Ricci curvature bounded below for Alexandrov spaces introduced by Zhang-Zhu,
we extend Theorem 1.1 to the case of Alexandrov spaces with an additional assumption.
\\

\noindent
{\bf Theorem 1.2}
\textit{
Given an $\epsilon>0$, an integer $n\geq 2$,
a non-increasing natural number valued function ${\mathcal{N}}(\ast)$,
and a positive function $\Phi_{0}(\star; \ast, n)$ depending on $\star$ and parameters $\ast, n$
with $\lim_{\star \rightarrow 0^{+}}\Phi_{0}(\star ; \ast, n)=0$,
there exists a $\delta=\delta(\epsilon, n, {\mathcal{N}}, \Phi_{0})>0$ such that if $M$ is an $n-$dimensional Alexandrov space
with $Ric(M)\geq n-1$ and $Vol(M)\geq\omega_{n}-\delta$,
and in addition,
$$
(A) \ \ \ \ \ \ \ \ \ \ \ \ \ \ \
\beta_{C_{p, \hat{\theta}}}(a)\leq\frac{{\mathcal{N}}(\hat{\theta})}{a^{n-2}}+\frac{\Phi_{0}[\omega_{n}-Vol(M); \hat{\theta}, n]}{a^{n-1}},
\ \ \ \ \ \  for \ any \ p\in M, \ \hat{\theta}, a\rightarrow 0^{+},
$$
then the Gromov-Hausdorff distance between $M$ and $S^{n}$ is at most $\epsilon$.
}
\\

\noindent
{\bf About the additional assumption (A)}
Roughly speaking, the assumption (A) says that the $n-1$ dimensional rough volume of the set of "short" cut points is small;
it is proposed to obtain Lemma 5.1 and make the proof of Lemma 5.5 go through;
while the assumption (A) seems reasonable by Lemma 4.3(ii) which says that
the $n-1$ dimensional Hausdorff volume of the set of directions of "short" cut points is small,
we do not know whether Lemma 5.1 and 5.5 could be deduced directly from Lemma 4.3(ii).
\\

%\noindent
%{\bf About the proof} \ \ Below are some possible approaches to prove Theorem 1.2.

%(a) Since the maximal diameter theorem on Alexandrov spaces [ZZ] obviously implies the maximal volume theorem,
%Theorem 1.2 is obtained if (a modification of) the definition of Ricci curvature bounded below for Alexandrov spaces introduced by Zhang-Zhu
%is stable with respect to Gromov-Hausdorff convergence.
%
%(b) Since CD is stable with respect to Gromov-Hausdorff convergence
%and Gigli......, if the cone of $M^{n}$ is an $n-$dimentional Alexandrov space with $Ric(M)\geq n-1$ .......[BBI]
%
%(c)or directly modify .....under CD

\noindent
{\bf About the proof}
The proof follows the basic idea in the Riemannian case with necessary modifications
because of the only almost everywhere second differentiability of distance functions.
More precisely, for the case of Riemannian manifolds, two key lemmas among others in the proof of Theorem 1.1 are
an $L^{2}-$version of Toponogov triangle comparison theorem for positive Ricci curvature (see Lemma 1.4 and 1.15 of [C])
which is deduced from Bochner's formula and a smooth approximation,
and a "predictability" lemma (see Lemma 2.10 of [C]) based on the $L^{2}-$version of Toponogov triangle comparison;
for the case of Alexandrov spaces, because of the lack of smoothness,
to obtain some kind of $L^{2}-$version of Toponogov triangle comparison for positive Ricci curvature,
we show an inequality of Riccati type for distance function and its solution along geodesics from some fixed point $p$,
and integrate this Riccati inequality outside the cut locus of $p$,
and thus obtain an $L^{2}-$version of Toponogov triangle comparison outside the cut locus of $p$,
then, to obtain a similar "predictability" lemma based on this Toponogov triangle comparison,
we have to bypass the cut locus of $p$ and thus the additional assumption (A) is proposed;
after having extended the two key lemmas from the case of Riemannian manifolds to the case of Alexandrov spaces,
the rest of the proof is the same as in the Riemannian case.
\\

\noindent
{\bf Acknowledgements} We would like to thank Professors Xi-Ping Zhu and Hui-Chun Zhang for their help and comments.
\\

{\section {Preliminaries}}

In this section, we recall basic concepts and some facts about Alexandrov spaces (see [BBI, BGP, PP, Pet2] for details),
and the Hessian and Laplacian for the distance function and the lower Ricci curvature bound for Alexandrov spaces introduced by Zhang-Zhu [ZZ1, ZZ2].

A complete metric space $(X, |\cdot  \cdot|)$ is called to be {\it a geodesic space} if for any two points $p, q\in X$,
the distance $|pq|$ is realized as the length of a rectifiable curve connecting $p$ and $q$.
Such distance-realizing curves, parameterized by arc-length, are called {\it (minimal) geodesics}.

Roughly speaking, given $k\in R$, a geodesic space $X$ is called to be {\it an Alexandrov space with curvature bounded from below by $k$ locally}
(for short, we say X to be {\it an Alexandrov space}), if all sufficiently small triangles in $X$ are not thinner than the corresponding $k-$plane triangles.
See [BBI, BGP, PP] for several equivalent definitions of Alexandrov spaces, here we recall the one from [PP]:

Let $\phi$ be a continuous function on $(a, b)$, $t\in (a, b)$. $\phi^{\prime\prime}(t)\leq B$ means
$\phi(t+\tau)\leq \phi(t) + A\tau + B\tau^{2}/2 + o(\tau^{2})$ for some $A\in R$.
If $\nu$ is another continuous function on $(a, b)$, then $\phi^{\prime\prime}\leq \nu$ means
$\phi^{\prime\prime}(t)\leq \nu(t)$ for all $t\in (a, b)$. Similarly, one can define
$\phi^{\prime\prime}(t)\geq B$ and $\phi^{\prime\prime}\geq \nu$;

$X$ is called {\it an Alexandrov space with curvature bounded from below by $k$} ($k\in R$) if
for any geodesic $\sigma: [0,l]\rightarrow X$ and $x$ not belonging to $\sigma[0,l]$,
$\phi:=\rho_{k}\circ dist_{x}\circ\sigma$ satisfies the differential inequality
$$
\phi^{\prime\prime}\leq 1-k\phi,
$$
where
$$
\rho_{k}(\varsigma):=
\left\{
\begin{array}{rl}
& 1/k[1-\cos(\varsigma\sqrt{k})],\ \ \ \ \ \ \ \ \ \ \  if \ \ k>0,
\\[3mm]
& x^{2}/2,\ \ \ \ \ \ \ \ \ \ \ \ \ \ \ \ \ \ \ \ \ \ \ \ \ \ \ \  if \ \ k=0,
\\[3mm]
& 1/k[1-\cosh(\varsigma\sqrt{-k})],\ \ \ \ \ \ \  if \ \ k<0.
\end{array}
\right.
$$

Below denote

$M$ --- be an $n-$dimensional Alexandrov space with curvature bounded below by $k_{0}(k_{0}<0)$ and without boundary,

$\uparrow_{x}^{y}$ --- the direction at $x$ of a minimal geodesic from $x$ to $y$,

$\Uparrow_{x}^{y}$ --- the set of directions at $x$ of minimal geodesics from $x$ to $y$,

$\Sigma_{x}$ --- the space of directions at a point $x\in M$,

$T_{x}$ --- the tangent cone at $x$,

$\sigma_{\xi}(s) := \exp_{x}(s\xi)(s\geq 0)$ --- some (quasi-)geodesic from $x$ in a direction $\xi\in\Sigma_{x}$.
\\

Fixing a point $p\in M$, let
$$
\gamma:[0,l]\rightarrow M \ \ be \ a\ geodesic\  with \ \gamma(0)=p \ \ and \ \ 0<l\leq\pi,
$$
$$
f:=dist_{p}=|p \cdot|, \ \ \ \ \ \ \ \ \ \ \ \ \tilde{f}:=\cos f,
$$
and $d_{x}f$ be the differential of $f$ at a point $x\in M\setminus\{p\} $.
By [BBI,4.5.7; BGP,11.4; OS,3.5], the first variation formula for distance function still holds, that is, for any $\xi\in\Sigma_{x}$,
$$
d_{x}f(\xi)=-\cos|\Uparrow_{x}^{p},\xi|.
$$
The lower and upper Hessian of $f$ at $x$ are defined respectively by
$$
\underline{Hess}_{x}f(\xi):=\liminf_{s\rightarrow 0}\frac{f(\exp_{x}(s\xi))-f(\exp_{x}(0))-d_{x}f(\xi)\cdot s}{s^{2}/2},
\ \ \ \ \ \ \ \ \ \ for \ any \ \ \xi\in\Sigma_{x},
$$
$$
\overline{Hess}_{x}f(\xi):=\limsup_{s\rightarrow 0}\frac{f(\exp_{x}(s\xi))-f(\exp_{x}(0))-d_{x}f(\xi)\cdot s}{s^{2}/2},
\ \ \ \ \ \ \ \ \ \ for \ any \ \ \xi\in\Sigma_{x},
$$
and $\underline{Hess}^{\varepsilon}_{x}f$, $\overline{Hess}^{\varepsilon}_{x}f$ are defined by replacing $\liminf_{s\rightarrow 0}$, $\limsup_{s\rightarrow 0}$ \ with \
$\liminf_{\varepsilon_{j}\rightarrow 0}$, $\limsup_{\varepsilon_{j}\rightarrow 0}$, respectively, for a sequence $\varepsilon:=\{\varepsilon_{j}\}_{j=1}^{\infty}$ with $\varepsilon_{j}\rightarrow 0$.

By the lower Alexandrov curvature bound $k_{0}(k_{0}<0)$ of $M$,
$$
\underline{Hess}_{x}f(\xi)\leq \overline{Hess}_{x}f(\xi)\leq
\frac{\sqrt{-k_{0}}\cosh(\sqrt{-k_{0}}| px |)}{\sinh(\sqrt{-k_{0}}| px |)}.
$$
Furthermore, by the triangle inequality,
\\

\noindent
{\bf Lemma 2.1 [ZZ1, (3.1, 3.2)]}
\textit{
(i) For any other geodesic $\sigma: [l_{1}, l_{2}]\rightarrow M$
with $\sigma[l_{1}, l_{2}]\cap\gamma[0,l]\ni \sigma(s_{0})=\gamma(t_{0})$ for some $s_{0}\in [l_{1}, l_{2}]$ and $t_{0}\in(0, l)$,
$$
-\frac{\sqrt{-k_{0}}\cosh[\sqrt{-k_{0}}(l-t_{0})]}{\sinh[\sqrt{-k_{0}}(l-t_{0})]}
\leq(f\circ\sigma)^{\prime\prime}(s_{0})
\leq\frac{\sqrt{-k_{0}}\cosh(\sqrt{-k_{0}}t_{0})}{\sinh(\sqrt{-k_{0}}t_{0})};
$$
(ii) for any $t\in(0,l)$ and $\xi\in\Sigma_{\gamma(t)}$,
$$
-\frac{\sqrt{-k_{0}}\cosh[\sqrt{-k_{0}}(l-t)]}{\sinh[\sqrt{-k_{0}}(l-t)]}
\leq\underline{Hess}_{\gamma(t)}f(\xi)\leq \overline{Hess}_{\gamma(t)}f(\xi)
\leq\frac{\sqrt{-k_{0}}\cosh(\sqrt{-k_{0}}t)}{\sinh(\sqrt{-k_{0}}t)}.
$$
}

\noindent
{\bf Remark 1} Lemma 2.1 (i) and (ii) are the same essentially, we write down different formulations for convenience.
\\

By [BGP, 10.6], for almost all $x\in M$, the tangent cone $T_{x}$ is isometric to Euclidean $n-$space $R^{n}$.
And according to [Per3], given a semiconcave function $\psi: M\rightarrow R$, for almost all $x\in M$,
there is a bilinear form $Hess_{x}\psi$ on $T_{x}$ such that
$$
\psi(y)=\psi(x)+d_{x}\psi(\uparrow_{x}^{y})\cdot|xy|+\frac{1}{2}Hess_{x}\psi(\uparrow_{x}^{y}, \uparrow_{x}^{y})\cdot|xy|^{2}+o(|xy|^{2})
$$
for any $\uparrow_{x}^{y}\in \Sigma_{x}$, moreover, $Hess_{x}\psi$ can be calculated using standard formulas,
hereafter, such an $x$ is called {\it a regular point} of $\psi$.
\\

\noindent
{\bf Remark 2}  If $\gamma(t)(0\leq t \leq l)$ is a regular point of $f$ and $\tilde{f}$,
since $Hess_{\gamma(t)}f$ and $Hess_{\gamma(t)}\tilde{f}$ can be calculated using standard formulas,
$\gamma^{\prime}(t)\in \Sigma_{\gamma(t)}$ is a eigenvector of $Hess_{\gamma(t)}f$ with the eigenvalue $0$,
and correspondingly, $\gamma^{\prime}(t)\in \Sigma_{\gamma(t)}$ is a eigenvector of $Hess_{\gamma(t)}\tilde{f}$ with the eigenvalue $-\cos t$.
\\

According to [BGP, Section 7], the tangent cone $T_{\gamma(t)}$ at an interior point $\gamma(t)(0<t<l)$ of $\gamma$
can be split into a direct metric product, denote
$$
L_{\gamma(t)}=\{\eta\in T_{\gamma(t)}|\angle(\eta,\gamma^{+}(t))=\angle(\eta,\gamma^{-}(t))=\pi/2\},
$$
$$
\Lambda_{\gamma(t)}=\{\eta\in\Sigma_{\gamma(t)}|\angle(\eta,\gamma^{+}(t))=\angle(\eta,\gamma^{-}(t))=\pi/2\}.
$$

The lower and upper Laplacian of $f$, \ lower and upper semi-Laplacian of $\tilde{f}$
at $\gamma(t)(0<t<l)$ are defined respectively by
$$
\underline{\triangle}_{\gamma(t)}f:=(n-1)\cdot\oint_{\Lambda_{\gamma(t)}}\underline{Hess}_{\gamma(t)}f(\eta),
\ \ \ \ \ \ \
\overline{\triangle}_{\gamma(t)}f:=(n-1)\cdot\oint_{\Lambda_{\gamma(t)}}\overline{Hess}_{\gamma(t)}f(\eta),
$$
$$
\underline{\diamondsuit}_{\gamma(t)}\tilde{f}:=(n-1)\cdot\oint_{\Lambda_{\gamma(t)}}\underline{Hess}_{\gamma(t)}\tilde{f}(\eta),
\ \ \ \ \ \ \
\overline{\diamondsuit}_{\gamma(t)}\tilde{f}:=(n-1)\cdot\oint_{\Lambda_{\gamma(t)}}\overline{Hess}_{\gamma(t)}\tilde{f}(\eta),
$$
and $\underline{\triangle}^{\varepsilon}_{\gamma(t)}f$, $\overline{\triangle}^{\varepsilon}_{\gamma(t)}f$ are defined
by replacing $\underline{Hess}_{\gamma(t)}f$, $\overline{Hess}_{\gamma(t)}f$ with
$\underline{Hess}^{\varepsilon}_{\gamma(t)}f$, $\overline{Hess}^{\varepsilon}_{\gamma(t)}f$, respectively, for a sequence $\varepsilon:=\{\varepsilon_{j}\}_{j=1}^{\infty}$ with $\varepsilon_{j}\rightarrow 0$.

If $\gamma(t)(0< t <l)$ is a regular point of $f$ and $\tilde{f}$,
the Laplacian of $f$, \ semi-Laplacian and Laplacian of $\tilde{f}$ at $\gamma(t)(0<t<l)$ are defined respectively by
$$
\triangle_{\gamma(t)}f:=(n-1)\cdot\oint_{\Lambda_{\gamma(t)}}Hess_{\gamma(t)}f(\eta),
$$
$$
\diamondsuit_{\gamma(t)}\tilde{f}:=(n-1)\cdot\oint_{\Lambda_{\gamma(t)}}Hess_{\gamma(t)}\tilde{f}(\eta),
\ \ \ \ \ \ \ \ \
\triangle_{\gamma(t)}\tilde{f}:=n\cdot\oint_{\Sigma_{\gamma(t)}}Hess_{\gamma(t)}\tilde{f}(\xi).
$$
where $\oint_{\diamond}\bullet:=\frac{1}{vol(\diamond)}\int_{\diamond}\bullet$.
\\

\noindent
{\bf Lemma 2.2}
\textit{
(i) For any other geodesic $\sigma: [l_{1}, l_{2}]\rightarrow M$
with $\sigma[l_{1}, l_{2}]\cap\gamma[0,l]\ni \sigma(s_{0})=\gamma(t_{0})$ for some $s_{0}\in [l_{1}, l_{2}]$ and $t_{0}\in(0, l)$,
$$
-\frac{\sqrt{-k_{0}}\cosh(\sqrt{-k_{0}}t_{0})}{\sinh(\sqrt{-k_{0}}t_{0})}-1
\leq(\tilde{f}\circ\sigma)^{\prime\prime}(s_{0})
\leq\frac{\sqrt{-k_{0}}\cosh[\sqrt{-k_{0}}(l-t_{0})]}{\sinh[\sqrt{-k_{0}}(l-t_{0})]}+1;
$$
(ii) For any $t\in (0, l)$ and $\eta\in\Lambda_{\gamma(t)}$,
$$
\underline{Hess}_{\gamma(t)}\widetilde{f}(\eta)=-\sin t \cdot \overline{Hess}_{\gamma(t)}f(\eta),
\ \ \ \
\overline{Hess}_{\gamma(t)}\widetilde{f}(\eta)=-\sin t \cdot \underline{Hess}_{\gamma(t)}f(\eta),
$$
$$
\underline{\diamondsuit}_{\gamma(t)}\widetilde{f}=-\sin t \cdot \overline{\triangle}_{\gamma(t)}f,
\ \ \ \ \ \ \ \
\overline{\diamondsuit}_{\gamma(t)}\widetilde{f}=-\sin t \cdot \underline{\triangle}_{\gamma(t)}f.
$$
}

\noindent
{\bf Proof.}
For simplification, let Lemma 2.1 (i) reads that $c_{1}\leq (f\circ\sigma)^{\prime\prime}(s_{0})\leq c_{2}$ with $-\infty<c_{1}<0<c_{2}<+\infty$,
that is, there exist $a_{1}, a_{2}\in R$ such that
$$
\begin{array}[b]{ll}
a_{1}(s-s_{0})+\frac{1}{2}c_{1}(s-s_{0})^{2}+o((s-s_{0})^{2})
&\leq (f\circ\sigma)(s)-(f\circ\sigma)(s_{0})
\\ \\
&\leq a_{2}(s-s_{0})+\frac{1}{2}c_{2}(s-s_{0})^{2}+o((s-s_{0})^{2}),
\end{array}
$$
by the first variation formula for the distance function, $d_{\sigma(s_{0})}f(\sigma^{\pm}(s_{0}))$ exist, thus
$$
a_{1}\geq -d_{\sigma(s_{0})}f(\sigma^{-}(s_{0}))\geq a_{2}\geq d_{\sigma(s_{0})}f(\sigma^{+}(s_{0}))\geq a_{1},
$$
and $f\circ\sigma$ is differentiable at $s_{0}$,    $(f\circ\sigma)^{\prime}(s_{0})=-d_{\sigma(s_{0})}f(\sigma^{-}(s_{0}))=d_{\sigma(s_{0})}f(\sigma^{+}(s_{0}))$, and
$$
\begin{array}[b]{ll}
\cos(f\circ\sigma)(s)=
&\cos(f\circ\sigma)(s_{0})-\sin(f\circ\sigma)(s_{0})\cdot[(f\circ\sigma)(s)-(f\circ\sigma)(s_{0})]
\\ \\
&+\frac{1}{2}[-\cos(f\circ\sigma)(s_{0})]\cdot[(f\circ\sigma)^{\prime}(s_{0})]^{2}\cdot(s-s_{0})^{2}+o((s-s_{0})^{2}),
\end{array}
$$
note that $0<(f\circ\sigma)(s_{0})<\pi$ and $|(f\circ\sigma)^{\prime}(s_{0})|\leq 1$, then,
$$
\begin{array}[b]{ll}
&-\sin(f\circ\sigma)(s_{0})\cdot(f\circ\sigma)^{\prime}(s_{0})\cdot(s-s_{0})
+\frac{1}{2}(-c_{2}-1)\cdot(s-s_{0})^{2}+o((s-s_{0})^{2})
\\ \\
\leq& \cos(f\circ\sigma)(s)-\cos(f\circ\sigma)(s_{0})
\\ \\
\leq&  -\sin(f\circ\sigma)(s_{0})\cdot(f\circ\sigma)^{\prime}(s_{0})\cdot(s-s_{0})
+\frac{1}{2}(-c_{1}+1)\cdot(s-s_{0})^{2}+o((s-s_{0})^{2}),
\end{array}
$$
that is, $-c_{2}-1\leq(\tilde{f}\circ\sigma)^{\prime\prime}(s_{0})\leq-c_{1}+1$, (i) is obtained.

Note that by the first variation formula, $d_{\gamma(t)}f(\eta)=0$ for any $t\in (0,l)$ and $\eta\in\Lambda_{\gamma(t)}$,
then, by Lemma 2.1 (ii), the proof of (ii) is similar to and easier than that of (i).
\\

In [Pet1], Petrunin proved the following second variation formula of arc-length,
based on which, Zhang-Zhu [ZZ1] introduced a new definition of Ricci curvature bounded below for Alexandrov spaces.
\\

\noindent
{\bf Proposition 2.3 [Pet1]}
\textit{
Given a geodesic $\gamma\subset M$, any two points $q_{1}, q_{2}\in\gamma$,which are not end points, and any positive number sequence $\{\tilde{\varepsilon}_{j}\}_{j=1}^{\infty}$
with $\tilde{\varepsilon}_{j}\rightarrow 0$, there exists a subsequence $\{\varepsilon_{j}\}\subset\{\widetilde{\varepsilon}_{j}\}$
and an isometry $T: L_{q_{1}}\rightarrow L_{q_{2}}$ such that
$$
|\exp_{q_{1}}(\varepsilon_{j}u), \exp_{q_{2}}(\varepsilon_{j}Tv)|
\leq |q_{1}q_{2}|+\frac{|uv|^{2}}{2|q_{1}q_{2}|}\cdot\varepsilon_{j}^{2}
\ \ \ \ \ \ \ \ \ \ \ \ \ \ \ \ \ \ \ \ \ \ \ \ \ \ \ \ \ \ \ \ \ \
$$
$$
\ \ \ \ \ \ \ \ \ \ \ \ \ \ \ \ \ \ \ \ \ \ \ \ \ \ \ \ \ \ \ \ \ \
-\frac{k_{0}\cdot|q_{1}q_{2}|}{6}\cdot(|u|^{2}+|v|^{2}+<u,v>)\cdot\varepsilon_{j}^{2}+o(\varepsilon_{j}^{2})
$$
for any $u, v\in L_{q_{1}}.$
}
\\

Note also that for a $2-$dimensional Alexandrov space, Cao-Dai-Mei [CDM]
improved the second variation formula such that the above inequality holds for all $\{\varepsilon_{j}\}_{j=1}^{\infty}$
with $\varepsilon_{j}\rightarrow 0$.
\\

\noindent
{\bf Definition 2.4 [ZZ2, Definition 2.5] (Condition $(RC)$)}
\textit{
Let $\alpha:[-s,s]\rightarrow M$ be a geodesic and $\{g_{\alpha(t)}\}_{-s<t<s}$ be a family of functions on $\Lambda_{\alpha(t)}$
such that $g_{\alpha(t)}$ is continuous on $\Lambda_{\alpha(t)}$ for each $t\in(-s,s)$.
The family $\{g_{\alpha(t)}\}_{-s<t<s}$ is said to satisfy Condition $(RC)$ on $\alpha$
if for any two points $q_{1}, q_{2}\in\alpha$
and any sequence $\{\tilde{\varepsilon}_{j}\}_{j=1}^{\infty}$ with $\tilde{\varepsilon}_{j}\rightarrow 0$,
there exists an isometry $T: L_{q_{1}}\rightarrow L_{q_{2}}$
and a subsequence $\{\varepsilon_{j}\}\subset\{\tilde{\varepsilon}_{j}\}$
such that
$$
|\exp_{q_{1}}(\varepsilon_{j}l_{1}\eta), \exp_{q_{2}}(\varepsilon_{j}l_{2}T\eta)|
\leq |q_{1}q_{2}|+\frac{(l_{1}-l_{2})^{2}}{2|q_{1}q_{2}|}\cdot\varepsilon_{j}^{2}
\ \ \ \ \ \ \ \ \ \ \ \ \ \ \ \ \ \ \ \ \ \ \ \ \ \ \ \ \ \ \ \ \ \ \ \ \ \ \
$$
$$
\ \ \ \ \ \ \ \ \ \ \ \ \ \ \ \ \ \ \ \ \ \ \ \ \ \ \ \ \ \ \ \ \ \ \ \
-\frac{g_{q_{1}}(\eta)\cdot|q_{1}q_{2}|}{6}\cdot[(l_{1})^{2}+l_{1}\cdot l_{2}+(l_{2})^{2}]\cdot\varepsilon_{j}^{2}
+o(\varepsilon_{j}^{2})
$$
for any $l_{1}, l_{2}\geq 0$ and any $\eta\in \Lambda_{q_{1}}.$
}
\\

\noindent
{\bf Definition 2.5 [ZZ2, Definition 2.6]}
\textit{
Let $\gamma:[0,l]\rightarrow M$ be a geodesic, $M$ is said to have Ricci curvature bounded below by $K$ along $\gamma$, if for any $\epsilon>0$ and any $0<t_{0}<l$,
there exists $\tau=\tau(t_{0},\epsilon)>0$
and a family of continuous functions $\{g_{\gamma(t)}\}_{t_{0}-\tau<t<t_{0}+\tau}$ on $\Lambda_{\gamma(t)}$
such that the family satisfies Condition $(RC)$ on $\gamma|_{(t_{0}-\tau, t_{0}+\tau)}$ and
$$
(n-1)\cdot\oint_{\Lambda_{\gamma(t)}}g_{\gamma(t)}(\eta)\geq K-\epsilon,
\ \ \ \ \ \ \ \ \forall t\in (t_{0}-\tau, t_{0}+\tau);
$$
M is said to have Ricci curvature bounded below by $K$, denoted by $Ric(M)\geq K$,
if each point $x\in M$ has a neighborhood $U_{x}$
such that $M$ has Ricci curvature bounded below by $K$ along every geodesic $\gamma$ in $U_{x}$.
}
\\

\noindent
{\bf Lemma 2.6  [ZZ1, Section 5, Para.3]}
\textit{For an n-dimensional Alexandrov space $M$ with $Ric(M)\geq n-1$ and without boundary,
the curvature-dimension condition $CD(n, n-1)$ holds and
$diam(M)\leq\pi$.
}
\\

The following lemma is a discrete version of the propagation equation of the Hessian of $f$ along the geodesic.
\\

\noindent
{\bf Lemma 2.7 (see also [ZZ1, Lemma 3.2])}
\textit{
Let $\gamma:[0,l]\rightarrow M$ be a geodesic with $\gamma(0)=p$,
$0<t_{0}-\tau<t_{1}<t_{0}<t_{2}<t_{0}+\tau<l$ with that $f$ is regular at $\gamma(t_{1})$ and $\gamma(t_{2})$,
and $\{g_{\gamma(t)}\}_{t_{0}-\tau<t<t_{0}+\tau}$ be a family of continuous functions on $\Lambda_{\gamma(t)}$
which satisfies Condition $(RC)$ on $\gamma|_{(t_{0}-\tau, t_{0}+\tau)}$,
then for any sequence $\{\widetilde{\varepsilon}_{j}\}_{j=1}^{\infty}$ with $\widetilde{\varepsilon}_{j}\rightarrow 0$,
there exists isometries $T_{1}: \Lambda_{\gamma(t_{0})}\rightarrow \Lambda_{\gamma(t_{1})}$
and $T_{2}: \Lambda_{\gamma(t_{0})}\rightarrow \Lambda_{\gamma(t_{2})}$,
and a subsequence $\varepsilon:=\{\varepsilon_{j}\}\subset\{\widetilde{\varepsilon}_{j}\}$,
such that
$$
\ \ \ \ \ \ (i)\ \ \ Hess_{\gamma(t_{1})}f(T_{1}\eta)
\geq l^{2}\cdot \overline{Hess}^{\varepsilon}_{\gamma(t_{0})}f(\eta)
-\frac{(l-1)^{2}}{t_{0}-t_{1}}
+\frac{l^{2}+l+1}{3}\cdot(t_{0}-t_{1})\cdot g_{\gamma(t_{0})}(\eta),
$$
$$
and \ (ii)\ \ \ Hess_{\gamma(t_{2})}f(T_{2}\eta)
\leq l^{2}\cdot \underline{Hess}^{\varepsilon}_{\gamma(t_{0})}f(\eta)
+\frac{(l-1)^{2}}{t_{2}-t_{0}}
-\frac{l^{2}+l+1}{3}\cdot(t_{2}-t_{0})\cdot g_{\gamma(t_{0})}(\eta)
$$
for any $\eta\in\Lambda_{\gamma(t_{0})}$ and any $l>0$.
}
\\

\noindent
{\bf Proof.}
The proof is a modification of that of [ZZ1, Lemma 3.2].

By the definition of Condition $(RC)$,
for the points $\gamma(t_{0}), \gamma(t_{1})$ and the sequence $\{\widetilde{\varepsilon}_{j}\}_{j=1}^{\infty}$,
there exists an isometry $T_{1}: \Lambda_{\gamma(t_{0})}\rightarrow \Lambda_{\gamma(t_{1})}$
and a subsequence $\varepsilon^{1}:=\{\varepsilon^{1}_{j}\}\subset\{\widetilde{\varepsilon}_{j}\}$,
such that the associated inequality holds;
Once again,
for the points $\gamma(t_{0}), \gamma(t_{2})$ and the sequence $\varepsilon^{1}:=\{\varepsilon^{1}_{j}\}_{j=1}^{\infty}$,
there exists an isometry $T_{2}: \Lambda_{\gamma(t_{0})}\rightarrow \Lambda_{\gamma(t_{2})}$
and a subsequence $\varepsilon:=\{\varepsilon_{j}\}\subset\{\varepsilon^{1}_{j}\}$,
such that the associated inequality holds.

For any $\eta\in\Lambda_{\gamma(t_{0})}$, $d_{\gamma(t_{0})}f(\eta)=0$,
choose a subsequence $\{\varepsilon^{\prime}_{j}\}_{j=1}^{\infty}\subset\{\varepsilon_{j}\}_{j=1}^{\infty}$ such that
$$\overline{Hess}^{\varepsilon}_{\gamma(t_{0})}f(\eta)
=\lim_{\varepsilon^{\prime}_{j}\rightarrow 0}
\frac{f(\exp_{\gamma(t_{0})}(\varepsilon^{\prime}_{j}\eta))-f(\gamma(t_{0}))}{(\varepsilon^{\prime}_{j})^{2}/2},$$
then for any $l>0$, let $\hat{\varepsilon}_{j}:=\frac{\varepsilon^{\prime}_{j}}{l}$, since $f$ is regular at $\gamma(t_{1})$,
$$
\begin{array}[b]{ll}
&f(\exp_{\gamma(t_{0})}(l\hat{\varepsilon}_{j}\eta))-f(\exp_{\gamma(t_{1})}(\hat{\varepsilon}_{j}T_{1}\eta))
\\ \\
=&[f(\gamma(t_{0}))+\frac{(l\hat{\varepsilon}_{j})^{2}}{2}\cdot \overline{Hess}^{\varepsilon}_{\gamma(t_{0})}f(\eta)+o((\hat{\varepsilon}_{j})^{2})]
\\ \\
&\ \ \ \ \ \ \ -[f(\gamma(t_{1}))+\frac{(\hat{\varepsilon}_{j})^{2}}{2}\cdot Hess_{\gamma(t_{1})}f(T_{1}\eta)+o((\hat{\varepsilon}_{j})^{2})]
\\ \\
=&t_{0}-t_{1}+\frac{(\hat{\varepsilon}_{j})^{2}}{2}\cdot
[l^{2}\cdot\overline{Hess}^{\varepsilon}_{\gamma(t_{0})}f(\eta)-Hess_{\gamma(t_{1})}f(T_{1}\eta)]+o((\hat{\varepsilon}_{j})^{2}),
\end{array}
$$
on the other hand, by the inequality associated with Condition $(RC)$,
$$
\begin{array}[b]{ll}
&f(\exp_{\gamma(t_{0})}(l\hat{\varepsilon}_{j}\eta))-f(\exp_{\gamma(t_{1})}(\hat{\varepsilon}_{j}T_{1}\eta))
\leq|\exp_{\gamma(t_{0})}(l\hat{\varepsilon}_{j}\eta),\exp_{\gamma(t_{1})}(\hat{\varepsilon}_{j}T_{1}\eta)|
\\ \\
&=t_{0}-t_{1}+\frac{(l-1)^{2}}{2(t_{0}-t_{1})}\cdot(\hat{\varepsilon}_{j})^{2}
-\frac{g_{\gamma(t_{0})}(\eta)\cdot(t_{0}-t_{1})}{6}\cdot(l^{2}+l+1)\cdot(\hat{\varepsilon}_{j})^{2}
+o((\hat{\varepsilon}_{j})^{2}),
\end{array}
$$
combining the above estimates and letting $\hat{\varepsilon}_{j}\rightarrow 0^{+}$, (i) is obtained.

Similarly, for any $\eta\in\Lambda_{\gamma(t_{0})}$,
choose a subsequence $\{\varepsilon^{\prime\prime}_{j}\}_{j=1}^{\infty}\subset\{\varepsilon_{j}\}_{j=1}^{\infty}$ such that
$$\underline{Hess}^{\varepsilon}_{\gamma(t_{0})}f(\eta)
=\lim_{\varepsilon^{\prime\prime}_{j}\rightarrow 0}
\frac{f(\exp_{\gamma(t_{0})}(\varepsilon^{\prime\prime}_{j}\eta))-f(\gamma(t_{0}))}{(\varepsilon^{\prime\prime}_{j})^{2}/2},$$
then for any $l>0$, let $\check{\varepsilon}_{j}:=\frac{\varepsilon^{\prime\prime}_{j}}{l}$, since $f$ is regular at $\gamma(t_{2})$,
$$
\begin{array}[b]{ll}
&f(\exp_{\gamma(t_{2})}(\check{\varepsilon}_{j}T_{2}\eta))-f(\exp_{\gamma(t_{0})}(l\check{\varepsilon}_{j}\eta))
\\ \\
=&[f(\gamma(t_{2}))+\frac{(\check{\varepsilon}_{j})^{2}}{2}\cdot Hess_{\gamma(t_{2})}f(T_{2}\eta)+o((\check{\varepsilon}_{j})^{2})]
\\ \\
&\ \ \ \ \ \ \ -[f(\gamma(t_{0}))+\frac{(l\check{\varepsilon}_{j})^{2}}{2}\cdot \underline{Hess}^{\varepsilon}_{\gamma(t_{0})}f(\eta)+o((\check{\varepsilon}_{j})^{2})]
\\ \\
=&t_{2}-t_{0}+\frac{(\check{\varepsilon}_{j})^{2}}{2}\cdot
[Hess_{\gamma(t_{2})}f(T_{2}\eta)-l^{2}\underline{Hess}^{\varepsilon}_{\gamma(t_{0})}f(\eta)]+o((\check{\varepsilon}_{j})^{2}),
\end{array}
$$
on the other hand, by the inequality associated with Condition $(RC)$,
$$
\begin{array}[b]{ll}
&f(\exp_{\gamma(t_{2})}(\check{\varepsilon}_{j}T_{2}\eta))-f(\exp_{\gamma(t_{0})}(l\check{\varepsilon}_{j}\eta))
\leq|\exp_{\gamma(t_{2})}(\check{\varepsilon}_{j}T_{2}\eta),\exp_{\gamma(t_{0})}(l\check{\varepsilon}_{j}\eta)|
\\ \\
&=t_{2}-t_{0}+\frac{(l-1)^{2}}{2(t_{2}-t_{0})}\cdot(\check{\varepsilon}_{j})^{2}
-\frac{g_{\gamma(t_{0})}(\eta)\cdot(t_{2}-t_{0})}{6}\cdot(l^{2}+l+1)\cdot(\check{\varepsilon}_{j})^{2}
+o((\check{\varepsilon}_{j})^{2}),
\end{array}
$$
combining the above estimates and letting $\check{\varepsilon}_{j}\rightarrow 0^{+}$, (ii) is obtained.
\\

{\section {An inequality of Riccati type}}

\noindent
{\bf Convention}
In this section,
let $\gamma:[0,l]\rightarrow M$ be a geodesic with $\gamma(0)=p$ and $0<l\leq \pi$ such that $f$ is regular almost everywhere on $\gamma$, and
$$
I:=(0,l),\ \ \ \ \ \ \ \ \ \
D:=\{t\in I | f \ is \ regular \ at \ \gamma(t)\},
$$
and $M$ has Ricci curvature bounded below by $n-1$.
\\

\noindent
{\bf Lemma 3.1 (see also [Pet3, Prop.2.2; ZZ1, Prop.A.4])}
\textit{
(i) $\triangle_{\gamma(t)}f$ is decreasing for $t\in D$,
that is, for any $t_{1},t_{2}\in D$ with $t_{1}<t_{2}$,
$$
\triangle_{\gamma(t_{1})}f
> \triangle_{\gamma(t_{2})}f;
$$
(ii) for any $t_{0}\in D$,
$$
\begin{array}[b]{ll}
\limsup_{D\ni t\rightarrow t_{0}}\frac{\triangle_{\gamma(t)}f-\triangle_{\gamma(t_{0})}f}{t-t_{0}}
&\leq -(n-1)-(n-1)\cdot\oint_{\Lambda_{\gamma(t_{0})}}(Hess_{\gamma(t_{0})}f(\eta))^{2}
\\ \\
&\leq -(n-1)-\frac{(\triangle_{\gamma(t_{0})}f)^{2}}{n-1}.
\end{array}
$$
}

\noindent
{\bf Proof.}
First to show that,
for any $\epsilon>0$, any sequence $\{\widetilde{\varepsilon}_{j}\}_{j=1}^{\infty}$ with $\widetilde{\varepsilon}_{j}\rightarrow 0^{+}$,
and any $t_{0}\in I$, there exists $\tau:=\tau(t_{0},\epsilon)>0$ such that the following holds:
for any $t_{1}\in (t_{0}-\tau, t_{0})\cap D$ and $t_{2}\in (t_{0}, t_{0}+\tau)\cap D$,
there exists a subsequence $\varepsilon:=\{\varepsilon_{j}\}\subset\{\widetilde{\varepsilon}_{j}\}$ such that
$$
\begin{array}[b]{ll}
\ \ \ \ \ \ (a)\ \ \ \ \ \ \
\frac{\triangle_{\gamma(t_{1})}f-\overline{\triangle}^{\varepsilon}_{\gamma(t_{0})}f}{t_{1}-t_{0}}
&\leq -(n-1)-\frac{n-1}{{\mathcal{H}}^{n-2}(\Lambda_{\gamma(t_{0})})}\cdot\int_{\Lambda_{\gamma(t_{0})}}(\overline{Hess}^{\varepsilon}_{\gamma(t_{0})}f(\eta))^{2}
+\epsilon
\\ \\
&\leq -(n-1)-\frac{(\overline{\triangle}^{\varepsilon}_{\gamma(t_{0})}f)^{2}}{n-1}+\epsilon,
\end{array}
$$
$$
\begin{array}[b]{ll}
and \ (b)\ \ \ \ \ \ \
\frac{\triangle_{\gamma(t_{2})}f-\underline{\triangle}^{\varepsilon}_{\gamma(t_{0})}f}{t_{2}-t_{0}}
&\leq -(n-1)-\frac{n-1}{{\mathcal{H}}^{n-2}(\Lambda_{\gamma(t_{0})})}\cdot\int_{\Lambda_{\gamma(t_{0})}}(\underline{Hess}^{\varepsilon}_{\gamma(t_{0})}f(\eta))^{2}
+\epsilon
\\ \\
&\leq -(n-1)-\frac{(\underline{\triangle}^{\varepsilon}_{\gamma(t_{0})}f)^{2}}{n-1}+\epsilon,
\end{array}
$$
then (i), (ii) are deduced from (a) and (b).
\\

Proof of (a) and (b):
by the definition of Ricci curvature bounded below by $n-1$, for any $\epsilon>0$ and any $t_{0}\in I$,
there exists $\tilde{\tau}=\tilde{\tau}(t_{0},\epsilon)>0$
and a family of continuous functions $\{g_{\gamma(t)}\}_{t_{0}-\tilde{\tau}<t<t_{0}+\tilde{\tau}}$ on $\Lambda_{\gamma(t)}$
such that the family satisfies Condition $(RC)$ on $\gamma|_{(t_{0}-\tilde{\tau}, t_{0}+\tilde{\tau})}$ and
$$
(n-1)\cdot\oint_{\Lambda_{\gamma(t)}}g_{\gamma(t)}(\eta)\geq (n-1)-\frac{\epsilon}{2},
\ \ \ \ \ \ \ \ \forall t\in (t_{0}-\tilde{\tau}, t_{0}+\tilde{\tau}).
$$

By Lemma 2.7, for any sequence $\{\widetilde{\varepsilon}_{j}\}_{j=1}^{\infty}$ with $\widetilde{\varepsilon}_{j}\rightarrow 0^{+}$,
any $t_{1}\in (t_{0}-\tilde{\tau}, t_{0})\cap D$ and $t_{2}\in (t_{0}, t_{0}+\tilde{\tau})\cap D$,
there exists isometries $T_{1}: \Lambda_{\gamma(t_{0})}\rightarrow \Lambda_{\gamma(t_{1})}$
and $T_{2}: \Lambda_{\gamma(t_{0})}\rightarrow \Lambda_{\gamma(t_{2})}$,
and a subsequence $\varepsilon:=\{\varepsilon_{j}\}\subset\{\widetilde{\varepsilon}_{j}\}$,
such that Lemma 2.7 (i) and (ii) hold.

To prove (a), consider
$$
\begin{array}[b]{ll}
F(l)
&:=l^{2}\cdot\overline{Hess}^{\varepsilon}_{\gamma(t_{0})}f(\eta)
-\frac{(l-1)^{2}}{t_{0}-t_{1}}
+\frac{l^{2}+l+1}{3}\cdot(t_{0}-t_{1})\cdot g_{\gamma(t_{0})}(\eta)
\\ \\
&=[\overline{Hess}^{\varepsilon}_{\gamma(t_{0})}f(\eta)+\frac{1}{t_{1}-t_{0}}-\frac{1}{3}(t_{1}-t_{0})\cdot g_{\gamma(t_{0})}(\eta)]\cdot l^{2}
\\ \\
& \ \ \ -[\frac{2}{t_{1}-t_{0}}+\frac{1}{3}(t_{1}-t_{0})\cdot g_{\gamma(t_{0})}(\eta)]\cdot l
\\ \\
& \ \ \ +[\frac{1}{t_{1}-t_{0}}-\frac{1}{3}(t_{1}-t_{0})\cdot g_{\gamma(t_{0})}(\eta)]
\\ \\
&=:A\cdot l^{2}+B\cdot l+C,
\end{array}
$$
by Lemma 2.1, $\overline{Hess}^{\varepsilon}_{\gamma(t_{0})}f(\eta)$ is bounded uniformly with respect to $\eta$,
so $A<0, B>0, -\frac{B}{2A}>0$, as $t_{1}\rightarrow t^{-}_{0}$.
Thus
$$
F(l)=A(l+\frac{B}{2A})^{2}+C-\frac{B^{2}}{4A},
$$
taking $l_{0}=-\frac{B}{2A}$,
$$
F(l_{0})=C-\frac{B^{2}}{4A}.
$$
Denote $\tau:=t_{1}-t_{0}, G:=g_{\gamma(t_{0})}(\eta), H:=\overline{Hess}^{\varepsilon}_{\gamma(t_{0})}f(\eta)$,
straightforward calculation shows
$$
\begin{array}[b]{ll}
F(l_{0})-H
&=(\frac{1}{\tau}-\frac{1}{3}\tau G)
-\frac{(\frac{2}{\tau}+\frac{1}{3}\tau G)^{2}}{4(H+\frac{1}{\tau}-\frac{1}{3}\tau G)}-H
\\ \\
&=-\frac{1}{3}G\cdot\tau+\frac{(-\frac{2}{3}G-H^{2})\cdot\tau+\frac{1}{3}GH\cdot\tau^{2}-\frac{1}{36}G^{2}\cdot\tau^{3}}
{1+H\cdot\tau-\frac{1}{3}G\cdot\tau^{2}}
\\ \\
&=-\frac{1}{3}G\cdot\tau-(\frac{2}{3}G+H^{2})\cdot\tau+\frac{(GH+H^{3})\cdot\tau-(\frac{1}{4}G^{2}+\frac{1}{3}GH^{2})\cdot\tau^{2}}
{1+H\cdot\tau-\frac{1}{3}G\cdot\tau^{2}}\cdot\tau
\\ \\
&=:-(G+H^{2})\cdot\tau+\varrho(\tau;t_{0},\eta)\cdot\tau,
\end{array}
$$
where $\lim_{\tau\rightarrow 0^{+}}\varrho(\tau;t_{0},\eta)=0$ uniformly with respect to $\eta$,
since $G, H$ are bounded uniformly with respect to $\eta$.

By Lemma 2.7 (i), and noting $\tau:=t_{1}-t_{0}<0$,
$$
\begin{array}[b]{ll}
\frac{Hess_{\gamma(t_{1})}f(T_{1}\eta)-\overline{Hess}^{\varepsilon}_{\gamma(t_{0})}f(\eta)}{t_{1}-t_{0}}
&\leq \frac{F(l_{0})-H}{\tau}
\\ \\
&=-(G+H^{2})+\varrho(\tau;t_{0},\eta)
\\ \\
&=-g_{\gamma(t_{0})}(\eta)-(\overline{Hess}^{\varepsilon}_{\gamma(t_{0})}f(\eta))^{2}+\varrho(\tau;t_{0},\eta),
\end{array}
$$
integrating over $\Lambda_{\gamma(t_{0})}$ on both sides of the above inequality,
and noting the isometry $T_{1}: \Lambda_{\gamma(t_{0})}\rightarrow \Lambda_{\gamma(t_{1})}$,
$$
\begin{array}[b]{ll}
\frac{\triangle_{\gamma(t_{1})}f-\overline{\triangle}^{\varepsilon}_{\gamma(t_{0})}f}{t_{1}-t_{0}}
\leq& -\frac{n-1}{{\mathcal{H}}^{n-2}(\Lambda_{\gamma(t_{0})})}\cdot\int_{\Lambda_{\gamma(t_{0})}}g_{\gamma(t_{0})}(\eta)
\\ \\
&-\frac{n-1}{{\mathcal{H}}^{n-2}(\Lambda_{\gamma(t_{0})})}\cdot\int_{\Lambda_{\gamma(t_{0})}}(\overline{Hess}^{\varepsilon}_{\gamma(t_{0})}f(\eta))^{2}
\\ \\
&+\frac{n-1}{{\mathcal{H}}^{n-2}(\Lambda_{\gamma(t_{0})})}\cdot\int_{\Lambda_{\gamma(t_{0})}}\varrho(\tau;t_{0},\eta).
\end{array}
$$
By H\"{o}lder inequality,
$$
\begin{array}[b]{ll}
&\frac{n-1}{{\mathcal{H}}^{n-2}(\Lambda_{\gamma(t_{0})})}\cdot\int_{\Lambda_{\gamma(t_{0})}}(\overline{Hess}^{\varepsilon}_{\gamma(t_{0})}f(\eta))^{2}
\\ \\
\geq &\frac{n-1}{{\mathcal{H}}^{n-2}(\Lambda_{\gamma(t_{0})})}\cdot\frac{1}{{\mathcal{H}}^{n-2}(\Lambda_{\gamma(t_{0})})}
\cdot(\int_{\Lambda_{\gamma(t_{0})}}\overline{Hess}^{\varepsilon}_{\gamma(t_{0})}f(\eta))^{2}
\\ \\
=&\frac{1}{n-1}
\cdot[\frac{n-1}{{\mathcal{H}}^{n-2}(\Lambda_{\gamma(t_{0})})}\cdot\int_{\Lambda_{\gamma(t_{0})}}\overline{Hess}^{\varepsilon}_{\gamma(t_{0})}f(\eta)]^{2}
\\ \\
=&\frac{(\overline{\triangle}^{\varepsilon}_{\gamma(t_{0})}f)^{2}}{n-1}.
\end{array}
$$
Denote
$\widetilde{\varrho}(\tau;t_{0}):=\frac{n-1}{{\mathcal{H}}^{n-2}(\Lambda_{\gamma(t_{0})})}\cdot\int_{\Lambda_{\gamma(t_{0})}}\varrho(\tau;t_{0},\eta)$,
then $\lim_{\tau\rightarrow 0^{+}}\widetilde{\varrho}(\tau;t_{0})=0$.
In particular, for the given $\epsilon$, there exists $0<\widehat{\tau}:=\widehat{\tau}(t_{0}, \epsilon)<\widetilde{\tau}$ such that
$|\widetilde{\varrho}(\tau;t_{0})|<\frac{\epsilon}{2}$ for any $t_{1}\in (t_{0}-\widehat{\tau}, t_{0})\cap D$.

Combining the above estimates, (a) is obtained.
\\

To prove (b),
similarly, consider
$$
\begin{array}[b]{ll}
F(l)&:= l^{2}\cdot \underline{Hess}^{\varepsilon}_{\gamma(t_{0})}f(\eta)
+\frac{(l-1)^{2}}{t_{2}-t_{0}}
-\frac{l^{2}+l+1}{3}\cdot(t_{2}-t_{0})\cdot g_{\gamma(t_{0})}(\eta)
\\ \\
&=[\underline{Hess}^{\varepsilon}_{\gamma(t_{0})}f(\eta)+\frac{1}{t_{2}-t_{0}}-\frac{1}{3}(t_{2}-t_{0})\cdot g_{\gamma(t_{0})}(\eta)]\cdot l^{2}
\\ \\
& \ \ \ -[\frac{2}{t_{2}-t_{0}}+\frac{1}{3}(t_{2}-t_{0})\cdot g_{\gamma(t_{0})}(\eta)]\cdot l
\\ \\
& \ \ \ +[\frac{1}{t_{2}-t_{0}}-\frac{1}{3}(t_{2}-t_{0})\cdot g_{\gamma(t_{0})}(\eta)]
\\ \\
&=:A\cdot l^{2}+B\cdot l+C,
\end{array}
$$
compared with the proof of (a), now the difference is $A>0, B<0$ as $t_{2}\rightarrow t^{+}_{0}$, and $\tau:=t_{2}-t_{0}>0$.
By Lemma 2.7 (ii), still
$$
\frac{Hess_{\gamma(t_{2})}f(T_{2}\eta)-\underline{Hess}^{\varepsilon}_{\gamma(t_{0})}f(\eta)}{t_{2}-t_{0}}
\leq \frac{F(l_{0})-H}{\tau}
$$
the remains are the same as that of (a).
\\

Proof of (i):
by (a) and (b), for any sequence $\{\widetilde{\varepsilon}_{j}\}_{j=1}^{\infty}$ with $\widetilde{\varepsilon}_{j}\rightarrow 0^{+}$,
and any $t_{0}\in I$, there exists $\tau:=\tau(t_{0})$ such that the following holds:
for any $t_{1}\in (t_{0}-\tau, t_{0})\cap D$ and $t_{2}\in (t_{0}, t_{0}+\tau)\cap D$,
there exists a subsequence $\varepsilon:=\{\varepsilon_{j}\}\subset\{\widetilde{\varepsilon}_{j}\}$ such that
$$
\triangle_{\gamma(t_{1})}f>\overline{\triangle}^{\varepsilon}_{\gamma(t_{0})}f
\geq\underline{\triangle}^{\varepsilon}_{\gamma(t_{0})}f>\triangle_{\gamma(t_{2})}f,
$$
then (i) follows from the compactness argument and the denseness of $D\subset I$.
\\

Proof of (ii):
just take $t_{0}\in D$ and let $t_{1}\rightarrow t_{0}^{-}, \ t_{2}\rightarrow t_{0}^{+}$ in (a) and (b).
\\

\noindent
{\bf Lemma 3.2 (see also [ZZ1, Theorem 3.3])}
\textit{
For any $t\in D$,
$$
(n-1)\cdot \cot(\pi-l+t)=-(n-1)\cdot \cot(l-t)
\leq\triangle_{\gamma(t)}f
\leq(n-1)\cdot \cot t.
$$
}

\noindent
{\bf Sketch of the proof}
(the details are contained in the appendix).

Denote $w(t):=\frac{\triangle_{\gamma(t)}f}{n-1}$,  then Lemma 3.1 reads that

(i) $w(t)$ is decreasing in $t\in D$,
that is, for any $t_{1},t_{2}\in D$ with  $t_{1}<t_{2}$,
$$
w(t_{1})>w(t_{2}),
$$
and (ii) for any $t_{0}\in D$,
$$
\limsup_{D\ni t\rightarrow t_{0}}\frac{w(t)-w(t_{0})}{t-t_{0}}
\leq -1-(w(t_{0}))^{2}.
$$
In addition,  since $f$ is regular on $D$, by Lemma 2.1,
$$
w(t)\leq\frac{\sqrt{-k_{0}}\cosh(\sqrt{-k_{0}}t)}{\sinh(\sqrt{-k_{0}}t)}\sim \frac{1}{t},
\ \ \ \ \ \ \ \ \ \  D\ni t\rightarrow 0^{+}.
$$
Thus, by successive approximations, for any $t\in D$,
$$
w(t)\leq\cot t,  \ \ \ \ \ \ \ \ \ \ \ \ \triangle_{\gamma(t)}f\leq(n-1)\cdot \cot t.
$$
Finally, by the triangle inequality, for any $t\in D$,
$$
\triangle_{\gamma(t)}f \geq -(n-1)\cdot \cot(l-t).
$$

\noindent
{\bf Lemma 3.3}
\textit{
Let $n\geq 2, \ 0<\theta\ll 1$, then there exists a positive function $\varphi(\theta)$ with $\lim_{\theta\rightarrow 0^{+}}\varphi(\theta)=0$ such that
$$
\begin{array}[b]{ll}
0
&\leq \int_{\theta}^{l(1-\theta)}
[(n-1)\cdot\oint_{\Lambda_{\gamma(t)}}(Hess_{\gamma(t)}\widetilde{f}(\eta))^{2}d\eta
-\frac{(\diamondsuit_{\gamma(t)}\tilde{f})^{2}}{n-1}]dt
\\ \\
&\leq
\left\{
\begin{array}{rl}
& \Phi(\theta; n),\ \ \ \ \ \ \  for \ \ \pi-\varphi(\theta)\leq l\leq \pi,
\\[3mm]
& \Psi(\theta; n),\ \ \ \ \ \ \  for \ \ \theta\leq l\leq \pi-\varphi(\theta).
\end{array}
\right.
\end{array}
$$
where $\lim_{\theta\rightarrow 0^{+}}\Phi(\theta; n)=0$, $\lim_{\theta\rightarrow 0^{+}}\Psi(\theta; n)=+\infty$.
}
\\

\noindent
{\bf Proof.}
Multiplying both sides of the inequality in Lemma 3.1(ii) by $\sin^{2}t$ and integrating with respect to $t\in [\theta,l(1-\theta)]$,
and by Lemma 2.2,
$$
\begin{array}[b]{ll}
0&\leq \int_{\theta}^{l(1-\theta)}
[(n-1)\cdot\oint_{\Lambda_{\gamma(t)}}(Hess_{\gamma(t)}\widetilde{f}(\eta))^{2}d\eta
-\frac{(\diamondsuit_{\gamma(t)}\tilde{f})^{2}}{n-1}]dt
\\ \\
&\leq \int_{\theta}^{l(1-\theta)}
[-(n-1)-\frac{(\triangle_{\gamma(t)}f)^{2}}{n-1}]\cdot\sin^{2}t \ dt
\\ \\
&\ \ \ -\int_{\theta}^{l(1-\theta)}
\limsup_{D\ni t^{\prime}\rightarrow t^{+}}\frac{\triangle_{\gamma(t^{\prime})}f-\triangle_{\gamma(t)}f}{t^{\prime}-t}
\cdot \sin^{2}t  \ dt
\\ \\
&=: I_{2}(\theta, l)-I_{1}(\theta, l).
\end{array}
$$

Now to estimate $I_{1}(\theta, l)$, take $\theta_{0}, l(1-\theta_{1})\in D$ with $\theta_{0}, \theta_{1}\rightarrow\theta^{+}$.

In case of $l(1-\theta) >\frac{\pi}{2}$, for any $N=1,2,...$, take a partition of $[\theta_{0}, l(1-\theta_{1})]$:
$\theta_{0}= t_{0}<t_{1}<...<t_{i}<t_{i+1}<...<t_{N_{0}-1}<t_{N_{0}}\leq\frac{\pi}{2}<t_{N_{0}+1}<...<t_{N}= l(1-\theta_{1})$
with $\{t_{i}\}_{i=0}^{N}\subset D$ and $\lim_{N\rightarrow \infty}[\max_{0\leq i\leq N-1}(t_{i+1}-t_{i})]=0$.
By Lemma 3.1 (i), $\triangle_{\gamma(t)}f$ is decreasing in $t\in D$, thus,
one can show that for any $i=0,1,..,N-1$,
$$
\int_{t_{i}}^{t_{i+1}}
\limsup_{D\ni t^{\prime}\rightarrow t^{+}}\frac{\triangle_{\gamma(t^{\prime})}f-\triangle_{\gamma(t)}f}{t^{\prime}-t}dt
\geq\triangle_{\gamma(t_{i+1})}f-\triangle_{\gamma(t_{i})}f,
$$
furthermore, by the uniformly continuity of $\sin t$ and Lemma 3.2,
$$
\begin{array}[b]{ll}
 &\int_{\theta_{0}}^{l(1-\theta_{1})}
\limsup_{D\ni t^{\prime}\rightarrow t^{+}}\frac{\triangle_{\gamma(t^{\prime})}f-\triangle_{\gamma(t)}f}{t^{\prime}-t}
\cdot \sin^{2}t  \ dt
\\ \\
\geq & \limsup_{N\rightarrow \infty}
\sum_{i=0}^{N-1}(\triangle_{\gamma(t_{i+1})}f-\triangle_{\gamma(t_{i})}f)\cdot\sin^{2}t_{i}
\\ \\
=&\limsup_{N\rightarrow \infty}
\{
\sum_{i=0}^{N-1}[\sin^{2}t_{i+1}\cdot\triangle_{\gamma(t_{i+1})}f-\sin^{2}t_{i}\cdot\triangle_{\gamma(t_{i})}f]
\\ \\
&\ \ \ \ \ \ \ \ \ \ \ \ \ \ \ \ \ \ +\sum_{i=0}^{N-1}(\sin^{2}t_{i}-\sin^{2}t_{i+1})\cdot\triangle_{\gamma(t_{i+1})}f
\}
\\ \\
=&\limsup_{N\rightarrow \infty}
\{
[\sin^{2}(l(1-\theta_{1}))\cdot\triangle_{\gamma(l(1-\theta_{1}))}f-\sin^{2}\theta_{0}\cdot\triangle_{\gamma(\theta_{0})}f]
\\ \\
&\ \ \ \ \ \ \ \ \ \ \ \ \ \ \ \ \ \ +\sum_{i=0}^{N-1}[\sin(2t_{i+1})\cdot(t_{i}-t_{i+1})+o(t_{i}-t_{i+1})]\cdot\triangle_{\gamma(t_{i+1})}f
\}
\\ \\
=&\limsup_{N\rightarrow \infty}
\{
[\sin^{2}(l(1-\theta_{1}))\cdot\triangle_{\gamma(l(1-\theta_{1}))}f-\sin^{2}\theta_{0}\cdot\triangle_{\gamma(\theta_{0})}f]
\\ \\
&\ \ \ \ \ \ \ \ \ \ \ \ \ \ \ \ \ \ +\sum_{i=0}^{N_{0}-1}[\sin(2t_{i+1})\cdot(t_{i}-t_{i+1})+o(t_{i}-t_{i+1})]\cdot\triangle_{\gamma(t_{i+1})}f
\\ \\
&\ \ \ \ \ \ \ \ \ \ \ \ \ \ \ \ \ \ +\sum_{i=N_{0}}^{N-1}[\sin(2t_{i+1})\cdot(t_{i}-t_{i+1})+o(t_{i}-t_{i+1})]\cdot\triangle_{\gamma(t_{i+1})}f
\}
\\ \\
\geq & \limsup_{N\rightarrow \infty}
\{
[\sin^{2}(l(1-\theta_{1}))\cdot(n-1)\cot(\pi-\theta_{1})-\sin^{2}\theta_{0}\cdot(n-1)\cot\theta_{0}]
\\ \\
&\ \ \ \ \ \ \ \ \ \ \ \ \ \ \ \ \ \ +\sum_{i=0}^{N_{0}-1}[\sin(2t_{i+1})\cdot(t_{i}-t_{i+1})+o(t_{i}-t_{i+1})]\cdot(n-1)\cot t_{i+1}
\\ \\
&\ \ \ \ \ \ \ \ \ \ \ \ \ \ \ \ \ \ +\sum_{i=N_{0}}^{N-1}[\sin(2t_{i+1})\cdot(t_{i}-t_{i+1})+o(t_{i}-t_{i+1})]\cdot(n-1)\cot(\pi-l+t_{i+1})
\}
\\ \\
=&[-\sin^{2}(l(1-\theta_{1}))\cdot(n-1)\cot\theta_{1}-\sin^{2}\theta_{0}\cdot(n-1)\cot\theta_{0}]
\\ \\
& \ -\int_{\theta_{0}}^{\frac{\pi}{2}}\sin(2t)\cdot(n-1)\cot t \ dt
-\int_{\frac{\pi}{2}}^{l(1-\theta_{1})}\sin(2t)\cdot(n-1)\cot(\pi-l+t) \ dt,
\end{array}
$$
letting $\theta_{0}, \theta_{1}\rightarrow\theta^{+}$,
$$
\begin{array}[b]{ll}
I_{1}(\theta, l)
&=\int_{\theta}^{l(1-\theta)}
\limsup_{D\ni t^{\prime}\rightarrow t^{+}}\frac{\triangle_{\gamma(t^{\prime})}f-\triangle_{\gamma(t)}f}{t^{\prime}-t}
\cdot \sin^{2}t  \ dt
\\ \\
&\geq -(n-1)
[\sin^{2}(l(1-\theta))\cdot\cot\theta+\sin\theta\cdot\cos\theta
\\ \\
& \ \ \ \ \ \ \ \ \ \ \ \ \ \ \ \ +\int_{\theta}^{\frac{\pi}{2}}\sin(2t)\cdot\cot t \ dt
+\int_{\frac{\pi}{2}}^{l(1-\theta)}\sin(2t)\cdot\cot(\pi-l+t) \ dt],
\\ \\
&=:J_{1}(\theta, l).
\end{array}
$$

Similarly, in case of $l(1-\theta) \leq \frac{\pi}{2}$, for any $N=1,2,...$, take a partition of $[\theta_{0}, l(1-\theta_{1})]$:
$\theta_{0}=t_{0}<t_{1}<...<t_{i}<t_{i+1}<...<t_{N}=l(1-\theta_{1})$, then
$$
\begin{array}[b]{ll}
 &\int_{\theta_{0}}^{l(1-\theta_{1})}
\limsup_{D\ni t^{\prime}\rightarrow t^{+}}\frac{\triangle_{\gamma(t^{\prime})}f-\triangle_{\gamma(t)}f}{t^{\prime}-t}
\cdot \sin^{2}t  \ dt
\\ \\
\geq&\limsup_{N\rightarrow \infty}
\{
[\sin^{2}(l(1-\theta_{1}))\cdot\triangle_{\gamma(l(1-\theta_{1}))}f-\sin^{2}\theta_{0}\cdot\triangle_{\gamma(\theta_{0})}f]
\\ \\
&\ \ \ \ \ \ \ \ \ \ \ \ \ \ \ \ \ \ +\sum_{i=0}^{N-1}[\sin(2t_{i+1})\cdot(t_{i}-t_{i+1})+o(t_{i}-t_{i+1})]\cdot\triangle_{\gamma(t_{i+1})}f
\}
\\ \\
\geq & \limsup_{N\rightarrow \infty}
\{
[\sin^{2}(l(1-\theta_{1}))\cdot(n-1)\cot(\pi-\theta_{1})-\sin^{2}\theta_{0}\cdot(n-1)\cot\theta_{0}]
\\ \\
&\ \ \ \ \ \ \ \ \ \ \ \ \ \ \ \ \ \ +\sum_{i=0}^{N-1}[\sin(2t_{i+1})\cdot(t_{i}-t_{i+1})+o(t_{i}-t_{i+1})]\cdot(n-1)\cot t_{i+1}
\}
\\ \\
=&[-\sin^{2}(l(1-\theta_{1}))\cdot(n-1)\cot\theta_{1}-\sin^{2}\theta_{0}\cdot(n-1)\cot\theta_{0}]
\\ \\
& \ -\int_{\theta_{0}}^{l(1-\theta_{1})}\sin(2t)\cdot(n-1)\cot t \ dt,
\end{array}
$$
letting $\theta_{0}, \theta_{1}\rightarrow \theta^{+}$,
$$
\begin{array}[b]{ll}
I_{1}(\theta, l)
&=\int_{\theta}^{l(1-\theta)}
\limsup_{D\ni t^{\prime}\rightarrow t^{+}}\frac{\triangle_{\gamma(t^{\prime})}f-\triangle_{\gamma(t)}f}{t^{\prime}-t}
\cdot \sin^{2}t  \ dt
\\ \\
&\geq -(n-1)
[\sin^{2}(l(1-\theta))\cdot\cot\theta +\sin\theta\cdot\cos\theta
+\int_{\theta}^{l(1-\theta)}2\cos^{2}t \ dt],
\\ \\
&=:J_{1}(\theta, l).
\end{array}
$$

Next to estimate $I_{2}(\theta, l)$. Note that generally, $I_{2}(\theta, l)<0$.

In case of $l(1-\theta) >\frac{\pi}{2}$,
by Lemma 3.2,
$$
|\triangle_{\gamma(t)}f|\geq
\left\{
\begin{array}{rl}
& \cot(\pi-l+t)>0,\ \ \ \ for \ \ 0< t< l-\frac{\pi}{2},
\\[3mm]
& 0,\ \ \ \ \ \ \ \ \ \ \ \ \ \ \ \ \ \ \ \ \ \ \ \ \ for \ \ l-\frac{\pi}{2}\leq t\leq \frac{\pi}{2},
\\[3mm]
& |\cot t|>0,\ \ \ \ \ \ \ \ \ \ \ \ \ for \ \ \frac{\pi}{2}< t< l,
\end{array}
\right.
$$
$$
\begin{array}[b]{ll}
I_{2}(\theta, l)
&=\int_{\theta}^{l(1-\theta)}[-(n-1)-\frac{(\triangle_{\gamma(t)}f)^{2}}{n-1}]\cdot\sin^{2}t \ dt
\\ \\
&\leq -(n-1)\cdot[\int_{\theta}^{l(1-\theta)}\sin^{2}t \ dt
\\ \\
&\ \ \ \ \ \ \ \ \ \ \ \ \ \ \ \ \ \ \ +\int_{\theta}^{l-\frac{\pi}{2}}\cot^{2}(\pi-l+t)\cdot\sin^{2}t \ dt
+\int_{\frac{\pi}{2}}^{l(1-\theta)}\cot^{2}t\cdot\sin^{2}t \ dt]
\\ \\
&=:J_{2}(\theta, l).
\end{array}
$$

Finally to estimate $I_{2}(\theta, l)-I_{1}(\theta, l)$.

To consider $J_{1}(\theta, l)$ in case of $l(1-\theta) >\frac{\pi}{2}$,
denote $H(\theta, l):=\int_{\frac{\pi}{2}}^{l(1-\theta)}\sin(2t)\cdot\cot(\pi-l+t) \ dt$,
since
$$
\lim_{\ l\geq \pi-\theta, \ \theta\rightarrow 0^{+}}\sin^{2}(l(1-\theta))\cdot \cot \theta = 0,
\ \ \ \ \ \
\lim_{\theta \rightarrow 0^{+}}H(\theta, \pi)=H(0, \pi)
$$
and for any $\theta>0$,
$$
\lim_{l \rightarrow \pi^{-}}H(\theta, l)=H(\theta, \pi),
$$
there exists a positive function $\varphi_{1}(\theta)$ with $\lim_{\theta\rightarrow 0^{+}}\varphi_{1}(\theta)=0$ such that
for $\pi-\varphi_{1}(\theta)\leq l\leq \pi$,
$$
\lim_{\theta\rightarrow 0^{+}}| J_{1}(\theta, l)-[-(n-1)\int_{0}^{\pi}\sin(2t)\cdot\cot t \ dt]|
=\lim_{\theta\rightarrow 0^{+}}| J_{1}(\theta, l)-[-(n-1)\pi]|=0.
$$
Similarly, consider $J_{2}(\theta, l)$ in case of $l(1-\theta) >\frac{\pi}{2}$,
there exists a positive function $\varphi_{2}(\theta)$ with $\lim_{\theta\rightarrow 0^{+}}\varphi_{2}(\theta)=0$ such that
for $\pi-\varphi_{2}(\theta)\leq l\leq \pi$,
$$
\lim_{\theta\rightarrow 0^{+}}| J_{2}(\theta, l)-[-(n-1)\int_{0}^{\pi}(\sin^{2}t+\cos^{2}t) \ dt]|
=\lim_{\theta\rightarrow 0^{+}}| J_{2}(\theta, l)-[-(n-1)\pi]|=0.
$$
Then, take $\varphi(\theta):=\min\{\theta, \varphi_{1}(\theta), \varphi_{2}(\theta)\}$, for $\pi-\varphi(\theta)\leq l\leq \pi$,
$$
0\leq\lim_{\theta\rightarrow 0^{+}}[I_{2}(\theta, l)-I_{1}(\theta, l)]\leq \lim_{\theta\rightarrow 0^{+}}[J_{2}(\theta, l)-J_{1}(\theta, l)]=0.
$$
Besides, note that generally,
$$
0\leq I_{2}(\theta, l)-I_{1}(\theta, l)
\leq 0 -J_{1}(\theta, l)
\leq (n-1)[\cot\theta+\theta+\pi+\frac{\pi}{2}\cdot\cot(\frac{\pi}{2}\theta)].
$$
Thus, the proof is completed.
\\

\noindent
{\bf Lemma 3.4}
\textit{
Let $n\geq 2, \ 0<\theta\ll 1$, then there exists a positive function $\varphi(\theta)$ with $\lim_{\theta\rightarrow 0^{+}}\varphi(\theta)=0$ such that
$$
\begin{array}[b]{ll}
0
&\leq \int_{\theta}^{l(1-\theta)}
[\oint_{\Sigma_{\gamma(t)}}(Hess_{\gamma(t)}\widetilde{f}(\xi))^{2}d\xi
-(\frac{1}{n}\triangle_{\gamma(t)}\widetilde{f})^{2}]dt
\\ \\
&\leq
\left\{
\begin{array}{rl}
& \Phi(\theta; n),\ \ \ \ \ \ \  for \ \ \pi-\varphi(\theta)\leq l\leq \pi,
\\[3mm]
& \Psi(\theta; n),\ \ \ \ \ \ \  for \ \ \theta\leq l\leq \pi-\varphi(\theta).
\end{array}
\right.
\end{array}
$$
where $\lim_{\theta\rightarrow 0^{+}}\Phi(\theta; n)=0$, $\lim_{\theta\rightarrow 0^{+}}\Psi(\theta; n)=+\infty$.
}
\\

\noindent
{\bf Proof.}
For any $t\in (\theta, l(1-\theta))\cap D$, denote all the eigenvalues of $Hess_{\gamma(t)}\tilde{f}$ by $\lambda_{1}=-\cos t, \lambda_{i}, 2\leq i\leq n$
(see Remark 2 in section 2),
then
$$\triangle_{\gamma(t)}\widetilde{f}=\sum_{1\leq i\leq n}\lambda_{i}
=-\cos t + \sum_{2\leq i\leq n}\lambda_{i},
\ \ \ \ \ \ \ \
\diamondsuit_{\gamma(t)}\widetilde{f}=\sum_{2\leq i\leq n}\lambda_{i}.
$$
By [HX, p.277],
$$
\begin{array}[b]{ll}
&\oint_{\Sigma_{\gamma(t)}}(Hess_{\gamma(t)}\widetilde{f}(\xi))^{2}d\xi
\\ \\
=&\frac{1}{n(n+2)}\cdot [\sum_{1\leq i\leq n}3\lambda_{i}^{2}+2\sum_{1\leq i<j\leq n}\lambda_{i}\lambda_{j}]
\\ \\
=&\frac{1}{n(n+2)}\cdot [(3\cos^{2} t-2\cos t\cdot\sum_{2\leq i\leq n}\lambda_{i})
+(\sum_{2\leq i\leq n}3\lambda_{i}^{2}+2\sum_{2\leq i<j\leq n}\lambda_{i}\lambda_{j})],
\end{array}
$$

$$
\oint_{\Lambda_{\gamma(t)}}(Hess_{\gamma(t)}\widetilde{f}(\eta))^{2}d\eta
=\frac{1}{(n-1)(n+1)}\cdot [\sum_{2\leq i\leq n}3\lambda_{i}^{2}+2\sum_{2\leq i<j\leq n}\lambda_{i}\lambda_{j}].
$$
Thus,
$$
\begin{array}[b]{ll}
&\int_{\theta}^{l(1-\theta)}
[\oint_{\Sigma_{\gamma(t)}}(Hess_{\gamma(t)}\widetilde{f}(\xi))^{2}d\xi
-(\frac{1}{n}\triangle_{\gamma(t)}\widetilde{f})^{2}]dt
\\ \\
=&\int_{\theta}^{l(1-\theta)}
\{
\frac{1}{n(n+2)}\cdot[(3\cos^{2}t-2\cos t\cdot\diamondsuit_{\gamma(t)}\widetilde{f})
+(n-1)(n+1)\oint_{\Lambda_{\gamma(t)}}(Hess_{\gamma(t)}\widetilde{f}(\eta))^{2}d\eta]
\\ \\
&\ \ \ \ \ \ -[-\frac{1}{n}\cos t+\frac{1}{n}\diamondsuit_{\gamma(t)}\widetilde{f}]^{2}
\}dt
\\ \\
=&\frac{2(n-1)}{n^{2}(n+2)}\cdot \int_{\theta}^{l(1-\theta)}[\frac{1}{n-1}\diamondsuit_{\gamma(t)}\widetilde{f}+\cos t]^{2}dt
\\ \\
&+\frac{(n-1)(n+1)}{n(n+2)}\cdot \int_{\theta}^{l(1-\theta)}[\oint_{\Lambda_{\gamma(t)}}(Hess_{\gamma(t)}\widetilde{f}(\eta))^{2}d\eta
-(\frac{1}{n-1}\diamondsuit_{\gamma(t)}\widetilde{f})^{2}]dt
\\ \\
=:&I_{1}(\theta, l)+I_{2}(\theta; l),
\end{array}
$$

Now to estimate $I_{1}(\theta, l)$. By Lemma 2.2 and 3.2, for any $t\in (0,l)$,
$$
0=-\sin t\cdot \cot t +\cos t
\leq \frac{1}{n-1}\diamondsuit_{\gamma(t)}\widetilde{f}+\cos t
\leq -\sin t\cdot \cot(\pi-l+t)+\cos t,
$$
$$
I_{1}(\theta, l)\leq
\frac{2(n-1)}{n^{2}(n+2)}\cdot \int_{\theta}^{l(1-\theta)}
[-\sin t\cdot \cot(\pi-l+t)+\cos t]^{2}dt=:J(\theta, l),
$$
Since $\lim_{l\rightarrow \pi^{-}}J(\theta, l)=J(\theta, \pi)=0$ for any small $\theta>0$,
there exists a positive function $\varphi(\theta)$ with $\lim_{\theta\rightarrow 0^{+}}\varphi(\theta)=0$ such that
for $\pi-\varphi(\theta)\leq l\leq \pi$,
$$
\lim_{\theta\rightarrow 0^{+}}J(\theta, l)=0.
$$
Besides, note that generally, as $t\in [\theta, l(1-\theta)]$,
$$
|\sin t\cdot \cot(\pi-l+t)|\leq
\left\{
\begin{array}{rl}
& \cos(\pi-l+t)\leq 1,\ \ \ \ \ \ \ \ \ \ \ \ \ \ \ \ \ \ for \ \ \theta<\pi-l+t\leq\frac{\pi}{2},
\\[3mm]
& \sin t\cdot |\cot(\pi-l\theta)|\leq\cot(\theta^{2}),\ \ \ \ for \ \ \frac{\pi}{2}\leq \pi-l+t\leq \pi-l\theta,
\end{array}
\right.
$$
$$
J(\theta, l)\leq
\frac{2(n-1)}{n^{2}(n+2)}\cdot \pi[2+\cot(\theta^{2})]^{2}.
$$

Thus, by Lemma 3.3, the proof is completed.
\\

{\section {An $L^{2}$ version of Toponogov triangle comparison}}

\noindent
{\bf Convention}
Hereafter, let $Ric (M) \geq n-1$, $Vol(M)>\omega_{n}-\delta$.
\\

\noindent
{\bf Notation}
Let $p, x \in M$, $v\in \Sigma_{p}$, $\xi\in \Sigma_{x}$, $\Gamma\subset\Sigma_{x}$,  $\theta, \hat{\theta}\rightarrow 0^{+}$,
$0\leq r_{1}< r_{2}\leq \pi$, $0^{+}\leftarrow r<l\leq\pi$, $R>0$,
denote
$$
B_{x}(R):=\{y\in M\mid |xy|<R\},
\ \ \ \ \ \ \
B_{x}^{\Gamma}(R):=\{y\in M\mid |xy|<R, \ and \  \exists \ some \ \uparrow_{x}^{y}\in\Gamma\},
$$
$$
A_{x}[r_{1}, r_{2}]:= \{y\in M \mid r_{1}\leq |xy|\leq r_{2}\},
\
A_{x}^{\Gamma}[r_{1}, r_{2}]:= \{y\in M \mid r_{1}\leq |xy|\leq r_{2}, \exists \ some \ \uparrow_{x}^{y}\in \Gamma \},
$$

similarly, \ $B_{1}(R)$, $B_{1}^{\Gamma}(R)$, $A_{1}[r_{1}, r_{2}]$, $A_{1}^{\Gamma}[r_{1}, r_{2}]$ are defined for the $n-$sphere $S^{n}$ with sectional curvature one,
$$
c(v):=\sup\{t\geq 0 \mid |p \exp_{p}(tv)|=t\},
$$
(note that $c(\uparrow_{p}^{x})$ is the same for any choice of $\uparrow_{p}^{x}\in\Uparrow_{p}^{x}$, thus $c(\Uparrow_{p}^{x})$ is well defined.)
$$
W_{p, \theta}:=\{x\in M \mid |px|<(1-\theta)c(\Uparrow_{p}^{x})\}\setminus \bar{B}_{p}(\theta),
$$
$$
\Gamma_{p, \hat{\theta}}:=\{v\in \Sigma_{p}\mid 0<c(v)\leq \pi-\hat{\theta}\},
\ \ \ \ \ \ \
\Gamma_{p, \hat{\theta}}^{c}:=\{v\in \Sigma_{p}\mid c(v)> \pi-\hat{\theta}\},
$$
$$
C_{p, \hat{\theta}}:=\{x\in M \mid  |px|=c(\Uparrow_{p}^{x})\leq \pi-\hat{\theta}\},
$$
$$
V_{p, \theta, \hat{\theta}}:=\{x\in M \mid  |px|\geq(1-\theta)c(\Uparrow_{p}^{x}) \ and \  \Uparrow_{p}^{x}\subset\Gamma_{p, \hat{\theta}}\},
$$
$$
V_{p, \theta, \hat{\theta}}^{c}:=\{x\in M \mid  |px|\geq(1-\theta)c(\Uparrow_{p}^{x}) \ and \  \Uparrow_{p}^{x}\subset\Gamma_{p, \hat{\theta}}^{c}\},
$$
$$
c(\xi):=\sup\{t\geq 0 \mid |x \exp_{x}(t\xi)|=t \},
\ \ \ \ \ \ \ \
l_{\xi}:=\min\{l, c(\xi)\},
$$
$$
\Sigma_{x,p,\theta,r,l}:=\{\xi\in\Sigma_{x}\mid l_{\xi}>r \ and\  \sigma_{\xi}[r, l_{\xi}]\subset W_{p, \theta}\},
$$
$$
h_{\tilde{f}}(\xi, s):=(\tilde{f}\circ\sigma_{\xi})(r)\cos (s-r)
+\frac{(\tilde{f}\circ\sigma_{\xi})(l_{\xi})-(\tilde{f}\circ\sigma_{\xi})(r)\cos (l_{\xi}-r)}{\sin(l_{\xi}-r)}\sin (s-r),
\
s\in [r, l_{\xi}].
$$

$\Phi(\star; \ast; \bullet,\bullet,...)$ is a positive function depending on $\star, \ast$ and parameters $\bullet,\bullet,...$
with the property that given any $\epsilon>0$ and parameters $\bullet,\bullet,...$, there exists a $\hat{\ast}:=\ast(\epsilon, \bullet,\bullet,...)$
and a $\hat{\star}:=\star(\epsilon, \hat{\ast}, \bullet,\bullet,...)$ such that $\Phi(\hat{\star}; \hat{\ast}; \bullet,\bullet,...)<\epsilon$.

$\Phi(\star; \bullet,\bullet,...)$, $\Phi(\star_{1}, \star_{2}; \bullet,\bullet,...)$, $\Phi(\star; \ast_{1}, \ast_{2}; \bullet,\bullet,...)$ are similar.

$C(\star, \ast)$ is a constant depending on $\star, \ast$.
\\

\noindent
{\bf Lemma 4.1 [LV, O, S2]}
\textit{
Let $Ric (M) \geq n-1$,  $\Gamma\subset\Sigma_{p}=S^{n-1}$ be any measurable subset,
and $0<r<R\leq\pi$, $0<s<S\leq\pi$, $r\leq s$, $R\leq S$, then
}

\textit{
(i) (Bishop absolute volume comparison)
$$
{\mathcal{H}}^{n}B_{p}^{\Gamma}(r)\leq {\mathcal{H}}^{n}B_{1}^{\Gamma}(r)={\mathcal{H}}^{n-1}(\Gamma)\cdot\int_{0}^{r}\sin^{n-1}t \ dt,
$$
$$
{\mathcal{H}}^{n}A_{p}^{\Gamma}[r, R]\leq {\mathcal{H}}^{n}A_{1}^{\Gamma}[r, R]={\mathcal{H}}^{n-1}(\Gamma)\cdot\int_{r}^{R}\sin^{n-1}t \ dt;
$$
}

\textit{
(ii) (Bishop-Gromov relative volume comparison)
$$
\frac{{\mathcal{H}}^{n}B_{p}^{\Gamma}(r)}{{\mathcal{H}}^{n}B_{p}^{\Gamma}(R)}
\geq \frac{{\mathcal{H}}^{n}B_{1}^{\Gamma}(r)}{{\mathcal{H}}^{n}B_{1}^{\Gamma}(R)},
\ \ \ \ \ \ \ \
\frac{{\mathcal{H}}^{n}A_{p}^{\Gamma}[r, R]}{{\mathcal{H}}^{n}A_{p}^{\Gamma}[s, S]}
\geq \frac{{\mathcal{H}}^{n}A_{1}^{\Gamma}[r, R]}{{\mathcal{H}}^{n}A_{1}^{\Gamma}[s, S]}.
$$
}

%\noindent
%{\bf Proof.}  modify ?............
%\\

\noindent
{\bf Lemma 4.2 [C, p.185 (2.5, 2.6)]}

\textit{
(i)Suppose that $Vol(M)>\omega_{n}-\delta$, by Bishop absolute volume comparison, for all $p\in M$, there exists a $q\in M$ with $d(p,q)>\pi-\epsilon$
(with $\lim_{\delta\rightarrow 0}\epsilon(\delta, n)=0$);
}

\textit{
(ii)By Bishop-Gromov relative volume comparison, if $p, q \in M$ with $d(p,q)>\pi-\delta$, then
$\forall x\in M, d(p, x)+ d(x, q)- d(p, q)< \epsilon$ (with $\lim_{\delta\rightarrow 0}\epsilon(\delta, n)=0$).
}
\\

\noindent
{\bf Lemma 4.3}
\textit{
(i) $diam(V_{p, \theta, \hat{\theta}}^{c})\leq\Phi(\theta, \hat{\theta}; n)$
\ where $\lim_{\theta, \hat{\theta}\rightarrow 0^{+}}\Phi(\theta, \hat{\theta}; n)=0$;
}

\textit{
(ii) ${\mathcal{H}}^{n-1}\Gamma_{p, \hat{\theta}}\leq\Phi(\delta; \hat{\theta}, n)$
\ where $\lim_{\delta\rightarrow 0^{+}}\Phi(\delta; \hat{\theta}, n)=0$.
}
\\

\noindent
{\bf Proof.}
(i) For any $x\in V_{p, \theta, \hat{\theta}}^{c}$, there exists a geodesic $\overline{py}$ such that
$x\in \overline{py}$, $|py|=c(\Uparrow_{p}^{x})>\pi-\hat{\theta}$, $|px|\geq(1-\theta)c(\Uparrow_{p}^{x})>(1-\theta)(\pi-\hat{\theta})$,
and $|xy|\leq \theta c(\Uparrow_{p}^{x})\leq\theta\pi$.

For any $x_{1}\in V_{p, \theta, \hat{\theta}}^{c}$ with $x_{1}\neq x$, by Lemma 4.2 (ii),
$|px_{1}|+|x_{1}y|<|py|+\epsilon(\hat{\theta}, n)$,
$|x_{1}y|<|py|-|px_{1}|+\epsilon(\hat{\theta}, n)\leq\pi-(1-\theta)(\pi-\hat{\theta})+\epsilon(\hat{\theta}, n)$.

Thus, $|xx_{1}|\leq |xy|+|x_{1}y|\leq \hat{\theta}+\theta(2\pi-\hat{\theta})+\epsilon(\hat{\theta}, n)=:\Phi(\theta, \hat{\theta}; n)$.
\\

(ii) By Lemma 4.1(i),
$$
\begin{array}[b]{ll}
\omega_{n}-\delta<Vol(M)
&=Vol(B_{p}^{\Gamma_{p, \hat{\theta}}}(\pi-\hat{\theta}))+Vol(B_{p}^{\Gamma_{p, \hat{\theta}}^{c}}(\pi))
\\ \\
&\leq {\mathcal{H}}^{n-1}\Gamma_{p, \hat{\theta}}\cdot\int_{0}^{\pi-\hat{\theta}}\sin^{n-1}t \ dt
+{\mathcal{H}}^{n-1}\Gamma_{p, \hat{\theta}}^{c}\cdot\int_{0}^{\pi}\sin^{n-1}t \ dt
\\ \\
&\leq{\mathcal{H}}^{n-1}\Gamma_{p, \hat{\theta}}\cdot\int_{0}^{\pi-\hat{\theta}}\sin^{n-1}t \ dt
+(\omega_{n-1}-{\mathcal{H}}^{n-1}\Gamma_{p, \hat{\theta}})\cdot\int_{0}^{\pi}\sin^{n-1}t \ dt
\\ \\
&=\omega_{n}-{\mathcal{H}}^{n-1}\Gamma_{p, \hat{\theta}}\cdot\int_{\pi-\hat{\theta}}^{\pi}\sin^{n-1}t \ dt,
\end{array}
$$
thus,
${\mathcal{H}}^{n-1}\Gamma_{p, \hat{\theta}}\leq\frac{\delta}{\int_{\pi-\hat{\theta}}^{\pi}\sin^{n-1}t \ dt}=:\Phi(\delta; \hat{\theta}, n)$.
\\

\noindent
{\bf Lemma 4.4}
\textit{
$\int_{M}dx\cdot\int_{\Sigma_{x,p,\theta,r,l}}d\xi\cdot
\int_{r}^{l_{\xi}}[Hess_{\sigma_{\xi}(s)}\widetilde{f}(\sigma_{\xi}^{+}(s))+\cos f(\sigma_{\xi}(s))]^{2}ds
\leq\Phi(\delta; \theta; n)$.
}
\\

\noindent
{\bf Proof.}
Noting that $Ric(M)\geq n-1$ and $\tilde{f}$ is regular almost everywhere on $M$, by Lemmas 3.4, 4.1(i) and 4.3(ii),
$$
\begin{array}[b]{ll}
&\int_{W_{p, \theta}}[\oint_{\Sigma_{x}}(Hess_{x}\widetilde{f}(\xi))^{2}d\xi
-(\frac{1}{n}\triangle_{x}\widetilde{f})^{2}]dx
\\ \\
\leq&\int_{\Gamma_{p, \varphi(\theta)}\sqcup\Gamma_{p, \varphi(\theta)}^{c}}dv\cdot
\int_{\theta}^{c(v)\cdot(1-\theta)}
[\oint_{\Sigma_{\gamma_{v}(t)}}(Hess_{\gamma_{v}(t)}\widetilde{f}(\xi))^{2}d\xi
-(\frac{1}{n}\triangle_{\gamma_{v}(t)}\widetilde{f})^{2}]dt
\\ \\
\leq&[\Phi(\delta; \varphi(\theta), n)\cdot\Psi(\theta; n)+\omega_{n-1}\cdot\Phi(\theta; n)]
\\ \\
=:&\Phi(\delta; \theta; n).
\end{array}
$$
Furthermore, by Lemma 4.3(ii), and the expressions for $\triangle_{\gamma_{v}(t)}\widetilde{f}, \diamondsuit_{\gamma_{v}(t)}\widetilde{f}$
and the estimate for $I_{1}(\theta, l)$ in the proof of Lemma 3.4,
$$
\begin{array}[b]{ll}
&\int_{W_{p, \theta}}[\oint_{\Sigma_{x}}(Hess_{x}\widetilde{f}(\xi)+\cos f(x))^{2}d\xi]dx
\\ \\
=&\int_{W_{p, \theta}}[\oint_{\Sigma_{x}}(Hess_{x}\widetilde{f}(\xi))^{2}d\xi
-(\frac{1}{n}\triangle_{x}\widetilde{f})^{2}]dx
+\int_{W_{p, \theta}}[\frac{1}{n}\triangle_{x}\widetilde{f}+\cos f(x)]^{2}dx
\\ \\
\leq&\Phi(\delta; \theta; n)+
\int_{\Gamma_{p, \varphi(\theta)}\sqcup\Gamma_{p, \varphi(\theta)}^{c}}dv\cdot
\int_{\theta}^{c(v)\cdot(1-\theta)}
[\frac{1}{n}\triangle_{\gamma_{v}(t)}\widetilde{f}+\cos f(\gamma_{v}(t))]^{2}dt
\\ \\
=&\Phi(\delta; \theta; n)+
(\frac{n-1}{n})^{2}\int_{\Gamma_{p, \varphi(\theta)}\sqcup\Gamma_{p, \varphi(\theta)}^{c}}dv\cdot
\int_{\theta}^{c(v)\cdot(1-\theta)}
[\frac{1}{n-1}\diamondsuit_{\gamma_{v}(t)}\widetilde{f}+\cos t]^{2}dt
\\ \\
\leq&\Phi(\delta; \theta; n)+
(\frac{n-1}{n})^{2}[\Phi(\delta; \varphi(\theta), n)\cdot\Psi(\theta; n)
+\omega_{n-1}\cdot\Phi(\theta; n)]
\\ \\
=:&\Phi(\delta; \theta; n).
\end{array}
$$
Thus, since geodesics do not branch in an Alexandrov space with curvature bounded below,
$$
\begin{array}[b]{ll}
&\int_{M}dx\cdot\int_{\Sigma_{x,p,\theta,r,l}}d\xi\cdot
\int_{r}^{l_{\xi}}[Hess_{\sigma_{\xi}(s)}\widetilde{f}(\sigma_{\xi}^{+}(s))+\cos f(\sigma_{\xi}(s))]^{2}ds
\\ \\
\leq&l\int_{W_{p, \theta}}[\int_{\Sigma_{x}}(Hess_{x}\widetilde{f}(\xi)++\cos f(x))^{2}d\xi]dx
\\ \\
\leq&\frac{\pi}{\omega_{n-1}}\int_{W_{p, \theta}}[\oint_{\Sigma_{x}}(Hess_{x}\widetilde{f}(\xi)++\cos f(x))^{2}d\xi]dx
\\ \\
\leq&\frac{\pi}{\omega_{n-1}}\cdot\Phi(\delta; \theta; n)=:\Phi(\delta; \theta; n).
\end{array}
$$
\\

\noindent
{\bf Lemma 4.5}
\textit{
For any $x\in M$ and $\xi\in \Sigma_{x,p,\theta,r,l}$, $(\tilde{f}\circ\sigma_{\xi})(s)$ is differentiable and $(\tilde{f}\circ\sigma_{\xi})^{\prime}(s)$ is lipschitz for $s \in [r, l_{\xi}]$.
}
\\

\noindent
{\bf Proof.}
For any $s \in [r, l_{\xi}]$,
by Lemma 2.2 (i),
$$
\begin{array}[b]{ll}
-\frac{\sqrt{-k_{0}}\cosh[\sqrt{-k_{0}}(f\circ\sigma_{\xi})(s)]}{\sinh[\sqrt{-k_{0}}(f\circ \sigma_{\xi})(s)]}-1
\leq(\tilde{f}\circ\sigma_{\xi})^{\prime\prime}(s)
\leq\frac{\sqrt{-k_{0}}\cosh[\sqrt{-k_{0}}(c(\Uparrow_{p}^{\sigma_{\xi}(s)})-(f\circ\sigma_{\xi})(s))]}
{\sinh[\sqrt{-k_{0}}(c(\Uparrow_{p}^{\sigma_{\xi}(s)})-(f\circ\sigma_{\xi})(s))]}+1,
\end{array}
$$
and by the definitions of $\Sigma_{x,p,\theta,r,l}$ and $W_{p, \theta}$,
$$
(f\circ\sigma_{\xi})(s)\geq\theta,
\ \ \ \ \ \ \ \ \ \
c(\Uparrow_{p}^{\sigma_{\xi}(s)})-(f\circ\sigma_{\xi})(s)\geq\theta\cdot c(\Uparrow_{p}^{\sigma_{\xi}(s)})
\geq\theta\cdot (f\circ\sigma_{\xi})(s)\geq\theta\cdot\theta=\theta^{2},
$$
then,
$$
c_{1}(\theta, k_{0}):=-\frac{\sqrt{-k_{0}}\cosh(\sqrt{-k_{0}}\pi)}{\sinh(\sqrt{-k_{0}}\theta)}-1
\leq(\tilde{f}\circ\sigma_{\xi})^{\prime\prime}(s)
\leq\frac{\sqrt{-k_{0}}\cosh(\sqrt{-k_{0}}\pi)}
{\sinh(\sqrt{-k_{0}}\theta^{2})}+1=:c_{2}(\theta, k_{0}).
$$
By [PP, 1.3 (1)], $(\tilde{f}\circ\sigma_{\xi})(s)-\frac{1}{2}c_{2}s^{2}$ is concave.
Besides, from the proof of Lemma 2.2 (i), $(\tilde{f}\circ\sigma_{\xi})^{\prime}(s)$ exists. Then,
for any $r\leq s_{1}< s_{2}\leq l_{\xi}$,
$$
(\tilde{f}\circ\sigma_{\xi})^{\prime}(s_{2})-(\tilde{f}\circ\sigma_{\xi})^{\prime}(s_{1})\leq c_{2}(s_{2}-s_{1}),
$$
similarly, $(\tilde{f}\circ\sigma_{\xi})^{\prime}(s_{2})-(\tilde{f}\circ\sigma_{\xi})^{\prime}(s_{1})\geq c_{1}(s_{2}-s_{1})$.
Thus, $(\tilde{f}\circ\sigma_{\xi})^{\prime}(s)$ is lipschitz for $s \in [r, l_{\xi}]$.
\\

\noindent
{\bf Lemma 4.6 (see also [C, Lemma 1.4, 1.15])} (Toponogov triangle comparison)
\textit{
Let $\frac{\pi}{2}\leq l_{0}<\pi$, then for any $0<l\leq l_{0}$,
$$
(i)\ \ \ \ \ \ \ \ \ \ \
\int_{M}dx\cdot\int_{\Sigma_{x,p,\theta,r,l}}d\xi\cdot\int_{r}^{l_{\xi}}
[(\widetilde{f}\circ\sigma_{\xi})(s)-h_{\tilde{f}}(\xi, s)]^{2}ds
\leq\Phi(\delta; \theta; l_{0}, n),
$$
$$
(ii)\ \ \ \ \ \ \ \ \
\int_{M}dx\cdot\int_{\Sigma_{x,p,\theta,r,l}}d\xi\cdot\int_{r}^{l_{\xi}}
[(\widetilde{f}\circ\sigma_{\xi})^{\prime}(s)-\frac{\partial h_{\tilde{f}}}{\partial s}(\xi, s)]^{2}ds
\leq\Phi(\delta; \theta; l_{0}, n).
$$
}
\\

\noindent
{\bf Proof.} The proof is a modification of the argument in the last two paragraphs of that of [C, Lemma 1.4].

For any $x\in M$, $\xi\in \Sigma_{x,p,\theta,r,l}$ and $s \in [r, l_{\xi}]$, set
$$
h(s):=(\tilde{f}\circ\sigma_{\xi})(s)-h_{\tilde{f}}(\xi, s).
$$
Easily $h(r)=h(l_{\xi})=0$, and by Lemma 4.5, $h(s)\cdot h^{\prime}(s)$ is lipschitz for $s \in [r, l_{\xi}]$, then
$$
\int_{r}^{l_{\xi}}[(h^{\prime})^{2}+h\cdot h^{\prime\prime}]ds=(h\cdot h^{\prime})(l_{\xi})-(h\cdot h^{\prime})(r)=0,
\ \ \ \ \ \ \
-\int_{r}^{l_{\xi}}h\cdot h^{\prime\prime}ds=\int_{r}^{l_{\xi}}(h^{\prime})^{2}ds.
$$
Besides, since $h^{\prime}(s)$ is lipschitz for $s \in [r, l_{\xi}]$, the Wirtinger's inequality holds,
$$
-\int_{r}^{l_{\xi}}h^{\prime\prime}\cdot hds=\int_{r}^{l_{\xi}}(h^{\prime})^{2}ds\geq(\frac{\pi}{l_{\xi}-r})^{2}\int_{r}^{l_{\xi}}h^{2}ds.
$$
Thus, by Lemma 4.4, using that $\frac{\partial^{2}h_{\tilde{f}}}{\partial s^{2}}=-h_{\tilde{f}}$ and integrating over $\Sigma_{x,p,\theta,r,l}$ and $M$ gives
$$
\begin{array}[b]{ll}
&(\frac{\pi}{l_{\xi}-r})^{2}\int_{M}dx\cdot\int_{\Sigma_{x,p,\theta,r,l}}d\xi\cdot\int_{r}^{l_{\xi}}
[(\widetilde{f}\circ\sigma_{\xi})(s)-h_{\tilde{f}}(\xi, s)]^{2}ds
\\ \\
\leq&-\int_{M}dx\cdot\int_{\Sigma_{x,p,\theta,r,l}}d\xi\cdot\int_{r}^{l_{\xi}}
[(\widetilde{f}\circ\sigma_{\xi})^{\prime\prime}(s)+h_{\tilde{f}}(\xi, s)][(\widetilde{f}\circ\sigma_{\xi})(s)-h_{\tilde{f}}(\xi, s)]ds
\\ \\
=&\int_{M}dx\cdot\int_{\Sigma_{x,p,\theta,r,l}}d\xi\cdot\int_{r}^{l_{\xi}}
[(\widetilde{f}\circ\sigma_{\xi})(s)-h_{\tilde{f}}(\xi, s)]^{2}ds
\\ \\
&-\int_{M}dx\cdot\int_{\Sigma_{x,p,\theta,r,l}}d\xi\cdot\int_{r}^{l_{\xi}}
[(\widetilde{f}\circ\sigma_{\xi})^{\prime\prime}(s)+(\widetilde{f}\circ\sigma_{\xi})(s)][(\widetilde{f}\circ\sigma_{\xi})(s)-h_{\tilde{f}}(\xi, s)]ds
\\ \\
\leq&\int_{M}dx\cdot\int_{\Sigma_{x,p,\theta,r,l}}d\xi\cdot\int_{r}^{l_{\xi}}
[(\widetilde{f}\circ\sigma_{\xi})(s)-h_{\tilde{f}}(\xi, s)]^{2}ds
\\ \\
&+[\Phi(\delta; \theta; n)]^{\frac{1}{2}}\cdot[\int_{M}dx\cdot\int_{\Sigma_{x,p,\theta,r,l}}d\xi\cdot\int_{r}^{l_{\xi}}
((\widetilde{f}\circ\sigma_{\xi})(s)-h_{\tilde{f}}(\xi, s))^{2}ds]^{\frac{1}{2}},
\end{array}
$$
note that $l_{\xi}-r<l\leq l_{0}<\pi$, one has
$$
\int_{M}dx\cdot\int_{\Sigma_{x,p,\theta,r,l}}d\xi\cdot\int_{r}^{l_{\xi}}
[(\widetilde{f}\circ\sigma_{\xi})(s)-h_{\tilde{f}}(\xi, s)]^{2}ds
\leq\frac{\Phi(\delta; \theta; n)}{[(\frac{\pi}{l_{0}})^{2}-1]^{2}}=:\Phi(\delta; \theta; l_{0}, n).
$$
(i) is obtained, contained in this set of inequalities is also (ii).
\\

\noindent
{\bf Corollary 4.7}
\textit{
Let $\frac{\pi}{2}\leq l_{0}<\pi, \ \{x_{i}\}_{i=1}^{m}\subset M, \ -1\leq\alpha_{i}\leq 1\ (i=1, 2, ...,m),
\ g:=\sum_{i=1}^{m}\alpha_{i}\cos dist_{x_{i}},\ h_{g}:=\sum_{i=1}^{m}\alpha_{i}h_{\cos dist_{x_{i}}}$,
 then for any $0<l\leq l_{0}$,
$$
(i)\ \ \ \ \ \ \ \ \ \ \
\int_{M}dx\cdot\int_{\cap_{i=1}^{m}\Sigma_{x,x_{i},\theta,r,l}}d\xi\cdot\int_{r}^{l_{\xi}}
[(g\circ\sigma_{\xi})(s)-h_{g}(\xi, s)]^{2}ds
\leq\Phi(\delta; \theta; l_{0}, m, n),
$$
$$
(ii)\ \ \ \ \ \ \ \ \
\int_{M}dx\cdot\int_{\cap_{i=1}^{m}\Sigma_{x,x_{i},\theta,r,l}}d\xi\cdot\int_{r}^{l_{\xi}}
[(g\circ\sigma_{\xi})^{\prime}(s)-\frac{\partial h_{g}}{\partial s}(\xi, s)]^{2}ds
\leq\Phi(\delta; \theta; l_{0}, m, n).
$$
}
\\

{\section {Proof of the theorem}}

\noindent
{\bf Notation and convention}
In this section, let ${\mathcal{N}}(\ast)$ be a non-increasing natural number valued function ,
and $\Phi_{0}(\star; \ast, n)$ be a positive function depending on $\star$ and parameters $\ast, n$
with $\lim_{\star \rightarrow 0^{+}}\Phi_{0}(\star ; \ast, n)=0$,
and suppose
$$
\beta_{C_{p, \hat{\theta}}}(a)\leq\frac{{\mathcal{N}}(\hat{\theta})}{a^{n-2}}+\frac{\Phi_{0}[\omega_{n}-Vol(M); \hat{\theta}, n]}{a^{n-1}},
\ \ \ \ \ \  for \ any \ p\in M, \ \hat{\theta}, a\rightarrow 0^{+}.
$$
And denote

$\Uparrow(y, \Omega):= \{\uparrow_{y}^{z}\in \Sigma_{y}\mid z\in\Omega\}$, \ for any $y\in M$, and $\Omega\subset M$;
\\

$\Gamma_{y_{1}, B_{y_{2}}(r), s}:=\{v\in\Sigma_{y_{1}}\mid v\in \Uparrow_{y_{1}}^{y} \ and\ \ |y_{1}y|=s \ for \ some\ \ y\in B_{y_{2}}(r)\}$,
\ \ \ \ \ \ \ \ \
for any
\\

\noindent
$y_{1}, y_{2}\in M$, $0<r\ll |y_{1}y_{2}|$, and $s>0$.
\\

\noindent
And a function $g: M\rightarrow R$ is called an {\it elementary function} generated by $\{x_{i}\}_{i=1}^{m}\subset M $ if
$g=\sum_{i=1}^{m}\alpha_{i}\cos dist_{x_{i}}$ for some $-1\leq\alpha_{i}\leq 1 \ (i=1, 2, ...,m)$.
\\

\noindent
{\bf Lemma 5.1}
\textit{
For any $y\in M$ and $r\in(0, \frac{\pi}{50})$,
$$
{\mathcal{H}}^{n-1}\Uparrow(y, (M\setminus W_{p, \theta}) \cap A_{y}[r, \frac{3}{4}\pi])
\leq \Phi(\delta;\theta; r, {\mathcal{N}}, \Phi_{0}, n).
$$
}

\noindent
{\bf Proof.}
First to show that for any $\hat{\theta}>0$, ${\mathcal{H}}^{n-1}\Uparrow(y, V_{p, \theta, \hat{\theta}}\cap A_{y}[r, \frac{3}{4}\pi])
\leq \Phi(\delta;\theta; \hat{\theta}, r, {\mathcal{N}}, \Phi_{0}, n)$.

Take $a=\pi\theta$, denote $\beta:=\beta_{C_{p, \hat{\theta}}}(a)$,
let $\{x_{i}\in C_{p, \hat{\theta}}\mid 1\leq i \leq\beta\}$
be a largest possible set of points in $C_{p, \hat{\theta}}$ that are at least $a$ pairwise distant from each other, then
$$
C_{p, \hat{\theta}}\subset\cup_{i=1}^{\beta}\bar{B}_{x_{i}}(2\pi\theta),
\ \ \ \ \ \
V_{p, \theta, \hat{\theta}}\subset\cup_{i=1}^{\beta}\bar{B}_{x_{i}}(3\pi\theta),
$$
and let $i=1, 2, ..., N_{y, r}(\leq \beta)$ be such that $\bar{B}_{x_{i}}(3\pi\theta)\cap A_{y}[r, \frac{3}{4}\pi]\neq\emptyset$.
Below suppose $0<\theta\ll r$, then for any $1\leq i\leq N_{y, r}$,
$$
\frac{r}{2}<r-6\pi\theta\leq|yx_{i}|-3\pi\theta\leq|yx_{i}|+3\pi\theta\leq\frac{3}{4}\pi+6\pi\theta<\frac{19}{25}\pi,
$$
denote $\Gamma_{i}:=\Uparrow(y, \bar{B}_{x_{i}}(3\pi\theta))$, then
$$
\bar{B}_{x_{i}}(3\pi\theta)\subset A_{y}^{\Gamma_{i}}[|yx_{i}|-3\pi\theta, |yx_{i}|+3\pi\theta]\subset \bar{B}_{x_{i}}(9\pi\theta),
$$
$$
\begin{array}[b]{ll}
\omega_{n}-\delta<Vol M&=Vol \bar{B}_{y}^{\Gamma_{i}}(|yx_{i}|-3\pi\theta)+Vol A_{y}^{\Gamma_{i}}[|yx_{i}|-3\pi\theta, |yx_{i}|+3\pi\theta]
\\ \\
& \ \ \ \
+Vol A_{y}^{\Gamma_{i}}[|yx_{i}|+3\pi\theta, \pi]
+Vol \bar{B}_{y}^{\Sigma_{y}\setminus\Gamma_{i}}(\pi)
\\ \\
&\leq{\mathcal{H}}^{n-1}\Gamma_{i}\cdot\int_{0}^{|yx_{i}|-3\pi\theta}\sin^{n-1}t dt+Vol \bar{B}_{x_{i}}(9\pi\theta)
\\ \\
& \ \ \ \
+{\mathcal{H}}^{n-1}\Gamma_{i}\cdot\int_{|yx_{i}|+3\pi\theta}^{\pi}\sin^{n-1}t dt
+(\omega_{n-1}-{\mathcal{H}}^{n-1}\Gamma_{i})\cdot\int_{0}^{\pi}\sin^{n-1}t dt
\\ \\
&\leq\omega_{n}-{\mathcal{H}}^{n-1}\Gamma_{i}\cdot\int_{|yx_{i}|-3\pi\theta}^{|yx_{i}|+3\pi\theta}\sin^{n-1}t dt
+\omega_{n-1}\cdot\int_{0}^{9\pi\theta}\sin^{n-1}t dt
\\ \\
&\leq\omega_{n}-{\mathcal{H}}^{n-1}\Gamma_{i}\cdot\sin^{n-1}\frac{r}{2}\cdot 6\pi\theta+\frac{\omega_{n-1}}{n}(9\pi\theta)^{n},
\end{array}
$$
thus,
${\mathcal{H}}^{n-1}\Gamma_{i}
\leq\frac{1}{6\sin^{n-1}\frac{r}{2}}\cdot(\frac{9^{n}\pi^{n-1}\omega_{n-1}}{n}\cdot\theta^{n-1}+\frac{1}{\pi}\cdot\frac{\delta}{\theta})$,
and
$$
\begin{array}[b]{ll}
&{\mathcal{H}}^{n-1}\Uparrow(y, V_{p, \theta, \hat{\theta}}\cap A_{y}[r, \frac{3}{4}\pi])
\leq\sum_{i=1}^{\beta}{\mathcal{H}}^{n-1}\Gamma_{i}
\\ \\
\leq& [\frac{{\mathcal{N}}(\hat{\theta})}{(\pi\theta)^{n-2}}+\frac{\Phi_{0}(\delta; \hat{\theta}, n)}{(\pi\theta)^{n-1}}]
\cdot \frac{1}{6\sin^{n-1}\frac{r}{2}}\cdot(\frac{9^{n}\pi^{n-1}\omega_{n-1}}{n}\cdot\theta^{n-1}+\frac{1}{\pi}\cdot\frac{\delta}{\theta})
\\ \\
=&\frac{1}{6\sin^{n-1}\frac{r}{2}}\cdot
[{\mathcal{N}}(\hat{\theta})\cdot(\frac{9^{n}\pi\omega_{n-1}}{n}\cdot\theta+\frac{1}{\pi^{n-1}}\cdot\frac{\delta}{\theta^{n-1}})
+\Phi_{0}(\delta; \hat{\theta}, n)\cdot(\frac{9^{n}\omega_{n-1}}{n}+\frac{1}{\pi^{n}}\cdot\frac{\delta}{\theta^{n}})]
\\ \\
=:&\Phi(\delta;\theta; \hat{\theta}, r, {\mathcal{N}}, \Phi_{0}, n).
\end{array}
$$

Next, similarly, if $\bar{B}_{p}(\theta)\cap A_{y}[r, \frac{3}{4}\pi]\neq\emptyset$,
$$
{\mathcal{H}}^{n-1}\Uparrow(y, \bar{B}_{p}(\theta)\cap A_{y}[r, \frac{3}{4}\pi])
\leq\frac{1}{3\sin^{n-1}\frac{r}{2}}\cdot(\frac{3^{n}\omega_{n-1}}{n}\cdot\theta^{n-1}+\frac{\delta}{\theta})
=:\Phi(\delta;\theta; r, n);
$$
and if $V_{p, \theta, \hat{\theta}}^{c}\cap A_{y}[r, \frac{3}{4}\pi]\neq\emptyset$, by Lemma 4.3 (i),
$V_{p, \theta, \hat{\theta}}^{c}\subset \bar{B}_{q}(\Phi(\theta, \hat{\theta}; n))$ for some $q\in M$,
$$
\begin{array}[b]{ll}
&{\mathcal{H}}^{n-1}\Uparrow(y, V_{p, \theta, \hat{\theta}}^{c}\cap A_{y}[r, \frac{3}{4}\pi])
\leq{\mathcal{H}}^{n-1}\Uparrow(y, \bar{B}_{q}(\Phi(\theta, \hat{\theta}; n))\cap A_{y}[r, \frac{3}{4}\pi])
\\ \\
\leq&\frac{1}{3\sin^{n-1}\frac{r}{2}}\cdot(\frac{3^{n}\omega_{n-1}}{n}\cdot\Phi^{n-1}(\theta, \hat{\theta}; n)
+\frac{\delta}{\Phi(\theta, \hat{\theta}; n)})
=:\Phi(\delta;\theta, \hat{\theta}; r, n).
\end{array}
$$

Finally, take $\hat{\theta}=\theta$, as
$M\setminus W_{p, \theta}=\bar{B}_{p}(\theta)\cup V_{p, \theta, \hat{\theta}}\cup V_{p, \theta, \hat{\theta}}^{c}$, Lemma 5.1 is obtained.
\\

\noindent
{\bf Lemma 5.2 (see also [C, Lemma 2.3])}
\textit{
For any $y_{1}, y_{2}\in M$, $0<r\ll |y_{1}y_{2}|$, and $s>0$, there exists some $s\in(|y_{1}y_{2}|-r, |y_{1}y_{2}|+r)$ such that
$$
{\mathcal{H}}^{n-1}\Gamma_{y_{1}, B_{y_{2}}(r), s}\geq\frac{nV_{n}(r)\cdot Vol M}{\pi^{n}\omega_{n}}=:C(r, n).
$$
where $V_{n}(r):= {\mathcal{H}}^{n}B_{1}(r)$.
}
\\

\noindent
{\bf Proof.} The proof is a modification of that of [C, Lemma 2.3].

Denote
$$
\Omega_{y_{1}, B_{y_{2}}(r)}:=\{v\in T_{y_{1}}\mid \frac{v}{|y_{1}y|}\in\Uparrow_{y_{1}}^{y} \ for \ some\ y\in B_{y_{2}}(r)\},
$$
by $Ric(M)\geq n-1$, $\exp_{y_{1}}: T_{y_{1}}\rightarrow M$ is volume non-increasing, and by Lemma 4.1 (ii),
$$
{\mathcal{H}}^{n}\Omega_{y_{1}, B_{y_{2}}(r)}\geq {\mathcal{H}}^{n}B_{y_{2}}(r)\geq\frac{V_{n}(r)\cdot Vol M}{\omega_{n}},
$$
but
$$
{\mathcal{H}}^{n}\Omega_{y_{1}, B_{y_{2}}(r)}
=\int_{0}^{\pi}{\mathcal{H}}^{n-1}\Gamma_{y_{1}, B_{y_{2}}(r), s}\cdot s^{n-1}ds
\leq \frac{\pi^{n}}{n}\cdot\max_{0\leq s\leq \pi}{\mathcal{H}}^{n-1}\Gamma_{y_{1}, B_{y_{2}}(r), s},
$$
Lemma 5.2 is obtained.
\\

\noindent
{\bf Lemma 5.3}
\textit{
For any $y_{1}, y_{2}\in M$ with $0<|y_{1}y_{2}|<\frac{5\pi}{8}$, $0<r\ll \min\{\frac{\pi}{16}, |y_{1}y_{2}|\}$, and $y\in B_{y_{1}}(\frac{r}{8})$,
$l:=|y_{1}y_{2}|-\frac{r}{4}$, the following holds:
$$
{\mathcal{H}}^{n-1}\Gamma_{y, B_{y_{2}}(r), l}\cap\Sigma_{y, p, \theta, \frac{r}{8}, l}
\geq C_{5.2}(r, n)-\Phi(\delta;\theta; r, {\mathcal{N}}, \Phi_{0}, n)
$$
where $C_{5.2}(r, n)$ refer to the $C(r, n)$ in Lemma 5.2; and for any $v\in\Gamma_{y, B_{y_{2}}(r), l}$,
$$
\exp_{y}(\frac{r}{8}v)\in B_{y_{1}}(r),\ \ \ \ \ \ \ \ \exp_{y}(lv)\in B_{y_{2}}(r).
$$
}

\noindent
{\bf Proof.} For any $y\in B_{y_{1}}(\frac{r}{8})$, by Lemma 5.2,
there exists some $s(y)\in(|y_{1}y_{2}|-\frac{r}{8}, |y_{1}y_{2}|+\frac{r}{8})$ such that
${\mathcal{H}}^{n-1}\Gamma_{y_{1}, B_{y_{2}}(r), s(y)}\geq C(r, n)$,
note that on the one hand, for any $z\in S_{y}(r)$,
$$
|y_{1}z|\leq|y_{1}y|+|yz|\leq \frac{r}{8}+\frac{r}{8}=\frac{r}{4}<r,
$$
on the other hand,
$$
s(y)\geq|yy_{2}|-\frac{r}{8}\geq|y_{1}y_{2}|-|y_{1}y|-\frac{r}{8}\geq|y_{1}y_{2}|-\frac{r}{8}-\frac{r}{8}=l,
$$
$$
s(y)\leq|yy_{2}|+\frac{r}{8}\leq|y_{1}y|+|y_{1}y_{2}|+\frac{r}{8}\leq\frac{r}{8}+|y_{1}y_{2}|+\frac{r}{8}=|y_{1}y_{2}|+\frac{r}{4}
<\frac{5\pi}{8}+\frac{\pi}{32}<\frac{3\pi}{4},
$$
$$
0<s(y)-l\leq(|y_{1}y_{2}|+\frac{r}{4})-(|y_{1}y_{2}|-\frac{r}{4})=\frac{r}{2},
$$
and for any $v\in\Gamma_{y, B_{y_{2}}(\frac{r}{8}), s(y)}$,
$$
|y_{2}\exp_{y}(lv)|\leq|y_{2}\exp_{y}[s(y)v]|+|\exp_{y}[s(y)v]\exp_{y}(lv)|
\leq\frac{r}{8}+(s(y)-l)\leq\frac{r}{8}+\frac{r}{2}=\frac{5r}{8}<r,
$$
thus, by Lemma 5.1,
$$
\begin{array}[b]{ll}
{\mathcal{H}}^{n-1}\Gamma_{y, B_{y_{2}}(r), l}\cap\Sigma_{y, p, \theta, \frac{r}{8}, l}
&\geq{\mathcal{H}}^{n-1}\Gamma_{y, B_{y_{2}}(\frac{r}{8}), s(y)}-{\mathcal{H}}^{n-1}\Uparrow(y, (M\setminus W_{p, \theta})\cap A_{y}[\frac{r}{8}, \frac{3}{4}\pi])
\\ \\
&\geq C_{5.2}(r, n)-\Phi(\delta;\theta; r, {\mathcal{N}}, \Phi_{0}, n).
\end{array}
$$

\noindent
{\bf Corollary 5.4}
\textit{
Lemma 5.3 still holds for any $\{x_{i}\}_{i=1}^{m}\subset M $ with $\Sigma_{y, p, \theta, \frac{r}{8}, l}$
replaced by $\cap_{i=1}^{m}\Sigma_{y,x_{i},\theta,\frac{r}{8},l}$,
and $\Phi(\delta;\theta; r, {\mathcal{N}}, \Phi_{0}, n)$ replaced by $\Phi(\delta;\theta; r, {\mathcal{N}}, \Phi_{0}, m, n)$.
}
\\

\noindent
{\bf Lemma 5.5 (see also [C, Lemma 2.10])}
\textit{
Given $\epsilon>0$, $l_{0}\in[\frac{\pi}{2}, \pi)$ and $m$ an integer, \
there exists \ $\mu(l_{0}, n)$ \ and  $\delta(\epsilon, m, l_{0}, {\mathcal{N}}, \Phi_{0}, k_{0}, n)$ such that the following holds:
If $Vol(M)>\omega_{n}-\delta$, $y_{1}, y_{2}\in M$ with $|y_{1}y_{2}|\leq \min\{\frac{5\pi}{8}, l_{0}\}$
and $f_{j}, g_{j}$ are at most $2m$ elementary fuctions generated by $\{x_{i}\}_{i=1}^{m}\subset M$
with $|f_{j}(y_{1})-g_{j}(y_{1})|<\mu\epsilon$ and $|f_{j}(y_{2})-g_{j}(y_{2})|<\mu\epsilon$,
then there exists $\bar{y}_{1}\in B_{y_{1}}(\mu\epsilon), \ \bar{y}_{2}\in B_{y_{2}}(\mu\epsilon)$
and a geodesic $\sigma_{\bar{y}_{1}\bar{y}_{2}}$ between them of length $\bar{l}=|\bar{y}_{1}\bar{y}_{2}|$
so that for all $s\in[0, \bar{l}]$,
$$
|f_{j}(\sigma_{\bar{y}_{1}\bar{y}_{2}}(s))-f_{j}(\bar{y}_{1})\cos s
-\frac{f_{j}(\bar{y}_{2})-f_{j}(\bar{y}_{1})\cos\bar{l}}{\sin\bar{l}}\sin s|<\frac{\epsilon}{n+1},
$$
$$
|f_{j}(\sigma_{\bar{y}_{1}\bar{y}_{2}}(s))-g_{j}(\sigma_{\bar{y}_{1}\bar{y}_{2}}(s))|<\epsilon.
$$
}

\noindent
{\bf Proof.}
The proof is a modification of that of [C, Lemma 2.10].

Set
$$\mu=\frac{1}{(n+1)(1+2\pi+\frac{2}{\sin l_{0}})},\ \ \ \ \ \  \  \ \ \ \ \ r=\mu\epsilon.$$
Suppose $0<\delta<\frac{1}{2}\omega_{n}$, then
$$
C_{5.2}(r, n)=\frac{nV_{n}(r)\cdot Vol M}{\pi^{n}\omega_{n}}\geq\frac{nV_{n}(r)}{2\pi^{n}},
\ \ \ \ \ \ \ \ \
{\mathcal{H}}^{n}B_{y_{1}}(\frac{r}{8})\geq\frac{V_{n}(\frac{r}{8})}{\omega_{n}}\cdot Vol M\geq\frac{1}{2}V_{n}(\frac{r}{8}).
$$
For $l:=|y_{1}y_{2}|-\frac{r}{4}<\min\{\frac{5\pi}{8}, l_{0}\}$ and any $y\in B_{y_{1}}(\frac{r}{8})$, denote
$$
\hat{\Sigma}_{y}:=\cap_{i=1}^{m}\Sigma_{y,x_{i},\theta,\frac{r}{8},l}\cap\Gamma_{y, B_{y_{2}}(r), l}.
$$
By Corollary 4.7, for any $e=f_{j}, \ g_{j}, \ j=1,2, ..., m$,
$$
\begin{array}[b]{ll}
&\int_{B_{y_{1}}(\frac{r}{8})}dy\cdot\int_{\hat{\Sigma}_{y}}d\xi
\cdot\int_{\frac{r}{8}}^{l}
[(e\circ\sigma_{\xi})^{\prime}(s)-\frac{\partial h_{e}}{\partial s}(\xi, s)]^{2}ds
\\ \\
\leq&\int_{M}dy\cdot\int_{\cap_{i=1}^{m}\Sigma_{y,x_{i},\theta,\frac{r}{8},l}}d\xi
\cdot\int_{\frac{r}{8}}^{l_{\xi}}
[(e\circ\sigma_{\xi})^{\prime}(s)-\frac{\partial h_{e}}{\partial s}(\xi, s)]^{2}ds
\leq\Phi_{4.7}(\delta; \theta; l_{0}, m, n),
\end{array}
$$
denote
$$
D_{y, e}:=\{\xi\in\hat{\Sigma}_{y}\mid
\int_{\frac{r}{8}}^{l}[(e\circ\sigma_{\xi})^{\prime}(s)-\frac{\partial h_{e}}{\partial s}(\xi, s)]^{2}ds
<\mu^{2}\epsilon^{2}\},
$$
$$
Y:=\{y\in B_{y_{1}}(\frac{r}{8})\mid{\mathcal{H}}^{n-1}\hat{\Sigma}_{y}\setminus D_{y, e}
\geq\frac{1}{2}{\mathcal{H}}^{n-1}\hat{\Sigma}_{y}\},
$$
then by Corollary 5.4,
$$
{\mathcal{H}}^{n}Y\leq\frac{\Phi_{4.7}(\delta; \theta; l_{0}, m, n)}{\mu^{2}\epsilon^{2}\cdot\frac{1}{2}{\mathcal{H}}^{n-1}\hat{\Sigma}_{y}}
\leq\frac{\Phi_{4.7}(\delta; \theta; l_{0}, m, n)}{\mu^{2}\epsilon^{2}\cdot\frac{1}{2}
[C_{5.2}(r, n)-\Phi_{5.4}(\delta;\theta; r, {\mathcal{N}}, \Phi_{0}, m, n)]},
$$
thus, note $C_{5.2}(r, n)\geq\frac{nV_{n}(r)}{2\pi^{n}}$ \ and \ ${\mathcal{H}}^{n}B_{y_{1}}(\frac{r}{8})\geq\frac{1}{2}V_{n}(\frac{r}{8})$,
as \ $\Phi_{4.7}(\delta; \theta; l_{0}, m, n)$ \ and

\noindent
$\Phi_{5.4}(\delta;\theta; r, {\mathcal{N}}, \Phi_{0}, n)$ are sufficiently small,
one has $\hat{Y}:=B_{y_{1}}(\frac{r}{8})\setminus Y$ such that
$$
{\mathcal{H}}^{n}\hat{Y}\geq\frac{1}{2}{\mathcal{H}}^{n}B_{y_{1}}(\frac{r}{8}),
$$
and for any $y\in\hat{Y}$,
$$
{\mathcal{H}}^{n-1}\hat{\Sigma}_{y}\setminus D_{y, e}<\frac{1}{2}{\mathcal{H}}^{n-1}\hat{\Sigma}_{y},
\ \ \ \ \
{\mathcal{H}}^{n-1}\hat{\Sigma}_{y}\cap D_{y, e}\geq\frac{1}{2}{\mathcal{H}}^{n-1}\hat{\Sigma}_{y}.
$$
After $2m$ steps of the above argument for all $e=f_{j}, \ g_{j}, \ j=1,2, ..., m$, one has that
as $\Phi_{4.7}(\delta; \theta; l_{0}, m, n)$ and $\Phi_{5.4}(\delta;\theta; r, {\mathcal{N}}, \Phi_{0}, n)$ are sufficiently small,
there exists some $\hat{\hat{Y}}\subset B_{y_{1}}(\frac{r}{8})$ such that
$$
{\mathcal{H}}^{n}\hat{\hat{Y}}\geq(\frac{1}{2})^{2m}{\mathcal{H}}^{n}B_{y_{1}}(\frac{r}{8})>0,
$$
and  for any $y\in\hat{\hat{Y}}$,
$$
{\mathcal{H}}^{n-1}\hat{\Sigma}_{y}\cap_{i=1}^{m}D_{y, f_{j}}\cap_{i=1}^{m}D_{y, g_{j}}\geq(\frac{1}{2})^{2m}{\mathcal{H}}^{n-1}\hat{\Sigma}_{y}>0,
$$
in particular, there exists some $\xi\in\hat{\Sigma}_{y}\cap_{i=1}^{m}D_{y, f_{j}}\cap_{i=1}^{m}D_{y, g_{j}}$ for some $y\in B_{y_{1}}(\frac{r}{8})$.
Let $\bar{y}_{1}:=\sigma_{\xi}(\frac{r}{8})\in B_{y_{1}}(r), \ \bar{y}_{2}:=\sigma_{\xi}(l)\in B_{y_{2}}(r)$,
then for all $e=f_{j}, \ g_{j}, \ j=1,2, ..., m$,
$$
\int_{\frac{r}{8}}^{l}|(e\circ\sigma_{\xi})^{\prime}(s)-\frac{\partial h_{e}}{\partial s}(\xi, s)|ds
\leq l^{\frac{1}{2}}\cdot[\int_{\frac{r}{8}}^{l}[(e\circ\sigma_{\xi})^{\prime}(s)-\frac{\partial h_{e}}{\partial s}(\xi, s)]^{2}ds]^{\frac{1}{2}}
<\pi^{\frac{1}{2}}\cdot\mu\epsilon,
$$
and for any $s\in[\frac{r}{8}, l]$,
$$
|e\circ\sigma_{\xi}(s)- h_{e}(\xi, s)|\leq |e\circ\sigma_{\xi}(\frac{r}{8})- h_{e}(\xi, \frac{r}{8})|+\pi^{\frac{1}{2}}\cdot\mu\epsilon
<\frac{\epsilon}{n+1},
$$
$$
|f_{j}\circ\sigma_{\xi}(s)-g_{j}\circ\sigma_{\xi}(s)|\leq|h_{f_{j}}(\xi, s)- h_{g_{j}}(\xi, s)|+2\pi^{\frac{1}{2}}\cdot\mu\epsilon
\leq\mu\epsilon+2\frac{\mu\epsilon}{\sin l_{0}}+2\pi^{\frac{1}{2}}\cdot\mu\epsilon<\epsilon.
$$
\\

\noindent
{\bf Lemma 5.6 (see also [Per4, Theorem 1])}
\textit{
For any integer $n\geq 2$ there exists $\delta=\delta(n)>0$ with the following property.
Let $M$ be an $n-$dimensional Alexandrov space without boundary and with $Ric\geq n-1$. Suppose $Vol(M)\geq \omega_{n}-\delta$.
Then $M$ is homeomorphic to $S^{n}$.
}
\\

\noindent
{\bf Proof.}
By [Per4], for the case of Riemannian manifolds, all the properties relevant to the proof
are a weakened version of the Abresch-Gromoll inequality
and a corollary of (the proof of) the Bishop-Gromov volume comparison inequality
(see [Per4, p.300] for details),
and the solution of (generalized) Poincar\'{e} conjecture [CZ, F, S].

For the case of Alexandrov spaces,
the Abresch-Gromoll inequality[GM, M, ZZ2] and the Bishop-Gromov volume comparison inequality [LV, O, S2] still hold,
the remain is to show that there exists $\delta=\delta(n)>0$ such that if $Vol(M)\geq\omega_{n}-\delta$
then $M$ is a topological manifold, which follows from

(a) By Lemma 4.1(i), given any $\epsilon>0$, there exists $\delta=\delta(\epsilon, n)>0$ such that if $Vol(M)\geq\omega_{n}-\delta$ then
$Vol(\Sigma_{x})\geq\omega_{n-1}-\epsilon$ for any $x\in M$;

(b) By [BBI, p.395, 10.9.15.3], given any $n\in N$ and $\epsilon>0$, there exists $\delta=\delta(\epsilon, n)>0$
such that for any $n-$dimensional Alexandrov space $X$ of curvature $\geq 1$
if $Vol(X)\geq\omega_{n}-\delta$ then $d_{GH}(X, S^{n})<\epsilon$;

(c) By [K, Per1], given any $n\in N$, there exists $\epsilon=\epsilon(n)>0$
such that for any $n-$dimensional Alexandrov space $X$ of curvature $\geq 1$
if $d_{GH}(X, S^{n})<\epsilon$ then $X$ is homeomorphic to $S^{n}$;

(d) By [Per1, Per2], for any $x\in M$, there exists a neighborhood of $x$ homeomorphic to
the tangent cone of $M$ at $x$.
\\

\noindent
{\bf Proof of the theorem} After having Lemma 5.5 and 5.6, the rest of the proof is the same as [C, 2.13-2.34].
\\

{\section {Appendix}}

\noindent
{\bf Proposition}
\textit{
Let $k_{0}<0$, $0<l\leq\pi$, $D\subset(0, l)$ be with full measure, $w: D\rightarrow R$ be a real valued function, satisfying
(i) $w(t)$ is decreasing in $t\in D$,
that is, for any $t_{1},t_{2}\in D$ such that  $t_{1}<t_{2}$,
$$
w(t_{1})>w(t_{2});
$$
(ii) for any $t_{0}\in D$,
$$
\limsup_{D\ni t\rightarrow t_{0}^{+}}\frac{w(t)-w(t_{0})}{t-t_{0}}
\leq -1-(w(t_{0}))^{2};
$$
and (iii)
$$
w(t)\leq\frac{\sqrt{-k_{0}}\cosh(\sqrt{-k_{0}}t)}{\sinh(\sqrt{-k_{0}}t)}\sim \frac{1}{t},
\ \ \ \ \ \ \ \ \ \  D\ni t\rightarrow 0^{+},
$$
then for any $t\in D$,
$$
w(t)\leq \cot t.
$$
}

\noindent
{\bf Proof.}
Since $w(t)$ is decreasing in $t\in D$, one may suppose that $w$ is differentiable almost everywhere on $D$. Denote
$$
E:=\{t\in D |\lim_{D\ni t\rightarrow t_{0}}\frac{w(t)-w(t_{0})}{t-t_{0}}=w^{\prime}(t_{0}) \ \  exists \ \ \},
$$
then $E\subset D\subset(0, l)$ is with full measure since $D\subset(0, l)$ is, in particular, $E$ is dense in $(0, l)$.
\\

Fixing $\bar{t}\in D$ and taking $D\ni t_{0}\rightarrow 0^{+}$ such that $w$ is differentiable at $t_{0}$, now to show that given any $\epsilon>0$,
$$
(\Box)\ \ \ \ \ \ \ \ \ \ \ \ \ \ \arctan w(\bar{t})-\arctan w(t_{0})\leq -(\bar{t}-t_{0})+2\epsilon\cdot(\bar{t}-t_{0}).
$$

By the differentiability of $w$ at $t_{0}$ and conditions (i) and (ii),
given any $\epsilon>0$, there exists $\tilde{\tau}_{0}:=\tilde{\tau}_{0}(t_{0}, \epsilon)$ such that for any $t\in (t_{0}, t_{0}+\tilde{\tau}_{0})\cap D$,
$$
\begin{array}[b]{ll}
\arctan w(t)-\arctan w(t_{0})
&<\frac{w(t)- w(t_{0})}{1+(w(t_{0}))^{2}}-\epsilon\cdot[w(t)- w(t_{0})]
\\ \\
&<\frac{w(t)- w(t_{0})}{1+(w(t_{0}))^{2}}-\epsilon\cdot[(w^{\prime}(t_{0})-\epsilon)\cdot(t-t_{0})],
\end{array}
$$
and
$$
\frac{w(t)-w(t_{0})}{t-t_{0}}< -[1+(w(t_{0}))^{2}]+\epsilon,
$$
$$
\frac{w(t)-w(t_{0})}{1+(w(t_{0}))^{2}}
<[-1+\frac{\epsilon}{1+(w(t_{0}))^{2}}]\cdot(t-t_{0})<(-1+\epsilon)\cdot(t-t_{0}),
$$
thus,
$$
\arctan w(t)-\arctan w(t_{0})<[(-1+\epsilon)-\epsilon\cdot(w^{\prime}(t_{0})-\epsilon)]\cdot(t-t_{0}).
$$
In short, given any $\epsilon>0$, there exists $\tau_{0}:=\tau_{0}(t_{0}, \epsilon)$ such that for any $t\in (t_{0}, t_{0}+\tau_{0})\cap D$,
$$
\arctan w(t)-\arctan w(t_{0})<-(t-t_{0})+\epsilon\cdot(t-t_{0}).
$$

Next take $t_{1}\in (t_{0}, t_{0}+\tau_{0})\cap D$ such that $w$ is differentiable at $t_{1}$, then
for the given $\epsilon>0$, there exists $\tau_{1}:=\tau_{1}(t_{1}, \epsilon)$ such that for any $t\in (t_{1}, t_{1}+\tau_{1})\cap D$,
$$
\arctan w(t)-\arctan w(t_{1})<-(t-t_{1})+\epsilon\cdot(t-t_{1}),
$$
and so on, one can take an increasing sequence $0<t_{0}<t_{1}<t_{2}<...<t_{n}<...<l$, such that for the given $\epsilon>0$,
$$
\begin{array}[b]{ll}
\arctan w(t_{1})-\arctan w(t_{0})
&<-(t_{1}-t_{0})+\epsilon\cdot(t_{1}-t_{0}),
\\ \\
\arctan w(t_{2})-\arctan w(t_{1})
&<-(t_{2}-t_{1})+\epsilon\cdot(t_{2}-t_{1}),
\\ \\
&......
\\ \\
\arctan w(t_{n})-\arctan w(t_{n-1})
&<-(t_{n}-t_{n-1})+\epsilon\cdot(t_{n}-t_{n-1}),
\end{array}
$$
summing up, for any $n=1, 2,...$,
$$
(\mho)\ \ \ \ \ \ \ \ \ \ \ \ \ \
\arctan w(t_{n})-\arctan w(t_{0})
<-(t_{n}-t_{0})+\epsilon\cdot(t_{n}-t_{0}).
$$

Below to show the inequality $(\Box)$ for three cases.

(a) In case of $\bar{t}\in (t_{n}, t_{n}+\tau_{n})\cap D$ for some $t_{n}$ taken above, just replace $t_{n+1}$ by $\bar{t}$,
the inequality $(\Box)$ is obtained by the above estimate $(\mho)$.

(b) In case of $\lim_{n\rightarrow \infty}t_{n}=\bar{t}$, by the estimate $(\mho)$ and condition (i), for any $n=1, 2,...$,
$$
\arctan w(\bar{t})-\arctan w(t_{0})
<\arctan w(t_{n})-\arctan w(t_{0})
<-(t_{n}-t_{0})+\epsilon\cdot(t_{n}-t_{0}),
$$
letting $n\rightarrow \infty$, the inequality $(\Box)$ is obtained.

(c) In case of $\lim_{n\rightarrow \infty}t_{n}=:t_{\infty}<\bar{t}$,
choose some $t_{n}$ such that $t_{\infty}-t_{n}$ is sufficiently small,
and take sufficiently small $\tau_{\infty}>0$ and $t_{*}\in(t_{\infty}, t_{\infty}+\tau_{\infty})$ such that $w$ is differentiable at $t_{*}$,
by condition (i),
$$
\arctan w(t_{*})-\arctan w(t_{n})<0=-(t_{*}-t_{n})+(t_{*}-t_{n}).
$$
Redefine $t_{n+1}:=t_{*}$, and denote $\bar{\epsilon}:=t_{*}-t_{n}$, then
$$
\arctan w(t_{n+1})-\arctan w(t_{n})<-(t_{n+1}-t_{n})+\bar{\epsilon},
$$
and since both $t_{\infty}-t_{n}$ and $\tau_{\infty}>0$ are sufficiently small,
one can suppose that $\bar{\epsilon}_{1}$ is arbitrarily small.
Starting again from $t_{n+1}:=t_{*}$, as above,
for the given $\epsilon>0$, there exists $\tau_{n+1}:=\tau_{n+1}(t_{n+1}, \epsilon)$ and $t_{n+2}\in (t_{n+1}, t_{n+1}+\tau_{n+1})\cap D$ such that
$$
\arctan w(t_{n+2})-\arctan w(t_{n+1})<-(t_{n+2}-t_{n+1})+\epsilon\cdot(t_{n+2}-t_{n+1}),
$$
and so on.
Finally, since $t_{n}-t_{n-1}>0$ for any $n\geq 1$, after at most countable steps,
one obtains an increasing sequence $0<t_{0}<t_{1}<t_{2}<...<t_{n}<...<l$ such that
$$
\lim_{n\rightarrow \infty}t_{n}=\bar{t},
$$
and for the given $\epsilon>0$ and any $n\geq 1$,
$$
\arctan w(t_{n})-\arctan w(t_{n-1})
<-(t_{n}-t_{n-1})+\epsilon\cdot(t_{n}-t_{n-1})+\bar{\epsilon}_{n},
$$
where $\bar{\epsilon}_{n}=(\frac{1}{2})^{n}\epsilon\cdot(\bar{t}-t_{0})$.
Summing up, for any $n=1, 2,...$,
$$
\arctan w(t_{n})-\arctan w(t_{0})
<-(t_{n}-t_{0})+\epsilon\cdot(t_{n}-t_{0})+\epsilon\cdot(\bar{t}-t_{0}),
$$
letting $n\rightarrow \infty$, the inequality $(\Box)$ is obtained.

Thus, the inequality $(\Box)$ is obtained for all cases.
\\

By the arbitrariness of $\epsilon$, one has that
for any $\bar{t}\in D$ and $t_{0}\in E$ with $t_{0}\rightarrow 0^{+}$,
$$
\arctan w(\bar{t})-\arctan w(t_{0})<-(\bar{t}-t_{0}),
$$
And by condition (iii), the proof is completed.
\\

\bigbreak

\noindent{\large\bf References}

\R{[AGS]} L. Ambrosio, N. Gigli, G. Savar\'{e}, Calculus and heat flow in metric measure spaces and applications to spaces with Ricci bounds from below, Invent. math., 195(2), 289¨C391 (2014).

\R{[BBI]} D. Burago, Y. Burago, S. Ivanov, A Course in Metric Geometry, Graduate Studies in Mathematics, vol. 33, AMS (2001).

\R{[BGP]} Y. Burago, M. Gromov, G. Perelman, A. D. Alexandrov spaces with curvatures bounded below, Russian Math. Surveys, 47, 1-58 (1992).

\R{[CDM]} J. G. Cao, B. Dai, J. Q. Mei, An optimal extension of Perelman's comparison theorem for quadrangles and its applications, Recent Advances in Geometric Analysis, ALM 11, 39-59 (2009).

\R{[CZ]} H. D. Cao, X. P. Zhu, Hamilton-Perelman's Proof of the Poincar\'{e} Conjecture and the Geometrization Conjecture, available at http://cn.arxiv.org/abs/0612069

\R{[C]} T. H. Colding, Shape of manifolds with positive Ricci curvature, Invent. math., 124, 175-191 (1996).

\R{[EKS]} M. Erbar, K. Kuwada, K-T. Sturm, On the Equivalence of the Entropic Curvature-Dimension Condition and Bochner's Inequality on Metric Measure Spaces,

http://cn.arxiv.org/abs/1303.4382.

\R{[F]} M. Freedman, The topology of four-manifolds, J.Differ. Geom., 17, 357-453 (1982).

\R{[G]} N. Gigli, The splitting theorem in non-smooth context, http://cn. arxiv. org/abs /1302. 5555

\R{[GM]} N. Gigli, S. Mosconi, The Abresch-Gromoll inequality in a non-smooth setting, http://cn. arxiv. org/abs/1209.3813

\R{[HX]} Z. S. Hu, S. L. Xu, An inequality between the integral norm and Euclidean norm of a symmetric bilinear form, J. of Math. Inequalities, 6 (2), 273-278 (2012).

\R{[K]}V. Kapovitch, Perelman's stability theorem, Surveys in differential geometry. Vol. XI, 103-136, Surv. Differ. Geom., 11, Int. Press, Somerville, MA, 2007.

\R{[KS1]} K. Kuwae, T. Shioya, Infinitesimal Bishop-Gromov condition for Alexandrov spaces, Probabilistic Approach to Geometry, 293¨C302, Adv. Stud. Prue Math. 57, Math. Soc. Japan, Tokyo, 2010.

\R{[KS2]} K. Kuwae, T. Shioya, A topological splitting theorem for weighted Alexandrov spaces, Tohoku Math. J., 63(2), no. 1, 59-76 (2011).

\R{[LV]} J. Lott, C. Villani, Ricci curvature for metric-measure spaces via optimal transport, Ann. of Math., 169, 903-991 (2009).

\R{[M]} M. Munn, Alexandrov spaces with large volume growth, http://cn.arxiv.org/abs /1405. 3312v1

\R{[O]} S. Ohta, On measure contraction property of metric measure spaces, Comment. Math. Helvetici, 82(4), 805-828 (2007).

\R{[Per1]} G. Perelman, A.D.Alexandrov¡¯s spaces with curvatures bounded from below, II, Preprint, available online at www.math.psu.edu/petrunin/

\R{[Per2]} G. Perelman, Elements of Morse theory on Alexandrov spaces, St. Petersburg Math. J., 5(1), 205¨C213 (1994).

\R{[Per3]} G. Perelman, DC structure on Alexandrov spaces. Preprint, preliminary version available online at www.math.psu.edu/petrunin/

\R{[Per4]} G. Perelman: Manifolds of positive Ricci curvature with almost maximal volume. JAMS 7, 299-305 (1994)

\R{[PP]} G. Perelman, A. Petrunin, Quasigeodesics and gradient curves in Alexandrov spaces, Preprint, available online at www.math.psu.edu/petrunin/

\R{[Pet1]} A. Petrunin, Parallel transportation for Alexandrov spaces with curvature bounded below, Geom. Funct. Analysis, 8(1), 123-148 (1998).

\R{[Pet2]} A. Petrunin, Semiconcave functions in Alexandrovs geometry, Surveys in Differential Geometry XI.

\R{[Pet3]} A. Petrunin, Alexandrov meets Lott-Villani-Sturm, M\"{u}nster J. of Math.,  4, 53-64 (2011).

\R{[QZZ]} Z. M. Qian, H. C. Zhang, X. P. Zhu, Sharp Spectral Gap and Li-Yau's Estimate on Alexandrov Spaces,  Mathematische Zeitschrift, 273 (3-4), 1175-1195 (2013).

\R{[S]} S. Smale, Generalized Poincar\'{e} Conjecture in dimensions greater than four, Ann. of Math., 74(2), 391-406 (1961).

\R{[S1]} K. Sturm, On the geometry of metric measure spaces. I.  Acta Math.,  196(1), 65-131 (2006).

\R{[S2]} K. Sturm, On the geometry of metric measure spaces. II.  Acta Math. 196(1), 133-177 (2006).

\R{[ZZ1]} H. C. Zhang, X. P. Zhu, Ricci curvature on Alexandrov spaces and rigidity theorems, Comm. Anal. Geom., 18(3), 503-554 (2010).

\R{[ZZ2]} H. C. Zhang, X. P. Zhu, On a new definition of Ricci curvature on Alexandrov spaces, Acta Mathematica Scientia, 30B(6), 1949-1974 (2010).

\R{[ZZ3]} H. C. Zhang, X. P. Zhu, Yau's gradient estimates on Alexandrov spaces, J. Differ. Geom., 91, 445-522 (2012).

\R{[ZZ4]} H. C. Zhang, X. P. Zhu, Lipschitz continuity of harmonic maps between Alexandrov spaces,  http://cn.arxiv.org/abs /1311.1331v3

\bigbreak

Zisheng Hu

Department of Mathematics and Computational Science

Shenzhen University

Shenzhen, Guangdong, 518060

China

email: zshu@szu.edu.cn

\bigbreak

Le Yin

Department of Mathematics and Computational Science

Shenzhen University

Shenzhen, Guangdong, 518060

China

email: lyin@szu.edu.cn

\bigbreak

\end{document}